\documentclass{amsart}
\usepackage{amssymb}
\ifx\pdfoutput\undefined
\usepackage{graphicx}  \DeclareGraphicsExtensions{.ps}

\else
\usepackage{amsfonts}
\usepackage[pdftex]{graphicx}  \DeclareGraphicsExtensions{.pdf,.mps}

\fi

\newcommand{\Op}{{\operatorname{Op}^{{w}}_h}}
\newcommand{\tR}{\widetilde R}
\newcommand{\tu}{\tilde u}
\newcommand{\tv}{\tilde v}

\newcommand{\CC}{{\mathbb C}}
\newcommand{\CI}{{\mathcal C}^\infty }
\newcommand{\CIc}{{\mathcal C}^\infty_{\rm{c}} }
 
\newcommand{\tH}{{\widetilde H}}

\newcommand{\ZZ}{{\mathbb Z}}

\newcommand{\RR}{{\mathbb R}}
\newcommand{\SP}{{\mathbb S}}
\newcommand{\DP}{{\mathbb D}}

\newcommand{\HH}{{\mathcal H}}

\newcommand{\C}{{\mathbb C}}
\newcommand{\NN}{{\mathbb N}}

\newcommand{\vol}{\operatorname{vol}}

\newcommand{\supp}{\operatorname{supp}}

\newcommand{\im}{\operatorname{Im}}
\newcommand{\defeq}{\stackrel{\rm{def}}{=}}

\newcommand{\Id}{\operatorname{Id}}
\newcommand{\tr}{\operatorname{tr}}
\newcommand{\dm}{\operatorname{dim}}
\newcommand{\ind}{\operatorname{ind}}

\newcommand{\rest}{\!\!\restriction}

\renewcommand{\Re}{\mathop{\rm Re}\nolimits}
\renewcommand{\Im}{\mathop{\rm Im}\nolimits}
\renewcommand{\im}{\mathop{\rm im}\nolimits}
\renewcommand{\ker}{\mathop{\rm ker}\nolimits}
\renewcommand{\dim}{\mathop{\rm dim}\nolimits}
\newcommand{\cker}{\mathop{\rm coker}\nolimits}

\theoremstyle{plain}

\newtheorem{thm}{Theorem}

\newtheorem{prop}{Proposition}[section]

\newtheorem{lem}[prop]{Lemma}

\theoremstyle{definition}

\numberwithin{equation}{section}

\def\bbbone{{\mathchoice {1\mskip-4mu {\rm{l}}} {1\mskip-4mu {\rm{l}}}
{ 1\mskip-4.5mu {\rm{l}}} { 1\mskip-5mu {\rm{l}}}}}

\def\squarebox#1{\hbox to #1{\hfill\vbox to #1{\vfill}}}

\title
{Elementary linear algebra for advanced spectral problems}
\author[J. Sj{\"o}strand]{Johannes Sj{\"o}strand}
\address{Centre de Math{\'e}matiques, {\'E}cole Polytechnique \\
UMR 7460, CNRS \\
F-91128 Palaiseau, France }
\email{johannes@math.polytechnique.fr}
\author[M. Zworski]{Maciej Zworski}
\address{Mathematics Department, University of California \\
Evans Hall, Berkeley, CA 94720, USA}
\email{zworski@math.berkeley.edu}

\begin{document}    
   
\maketitle   
   
\section{Introduction}   
\label{in}

The purpose of this article is to discuss a simple linear 
algebraic tool which has proved itself very useful in the
mathematical study of spectral problems arising in 
elecromagnetism and quantum mechanics. Roughly speaking 
it amounts to replacing an operator of interest by 
a suitably chosen invertible system of operators. 

That approach has a very long tradition and
appears constantly under different names and guises in 
many works of pure and applied mathematics. 
Our purpose here is not to provide a historical survey but
to present an account of a specific approach from a
personal perspective of the authors. On one hand we 
hope to provide a source of systematic references 
for the practitioners of our type of spectral theory
and, hopefully,  to convince others of the usefulness of this method.
We do not know, but find very interesting, if the method which has proved
itself so successful in theoretical studies has a chance of being useful
numerically.

The key elementary observation goes back --  at least -- to Schur and 
his {\em complement} formula: if for matrices
\[ \left( \begin{array}{ll} P & R_- \\
R_+ & 0 \end{array} \right)^{-1} = 
\left( \begin{array}{ll} E & E_+ \\E_- & E_{-+} \end{array}
\right) \,,\]
then $ P $ is invertible if and only if $ E_{-+} $ is invertible
and 
\begin{equation}
\label{eq:schur1} P^{-1} = E - E_+ E_{-+}^{-1} E_- \,, \ \ 
E_{ -+}^{-1}  = - R_- P^{-1} R_+  \,.\end{equation}
In fact the equivalence of invertibilities of $ P $ and
$ E_{-+} $ holds for systems with a non
zero lower right hand corner  (see Lemma \ref{l:schur}) 
but since here we always start with $ P$ and choose $ R_\pm $ we
can normally consider these simpler systems. Sometimes, in the context
of index theory one considers operators $ P $ which 
are never invertible. In that case the index of $ P $ is equal the index of
$ E_{-+ } $ which is trivial to compute  if $ E_{-+} $ is a matrix -- 
 see \S \ref{seft}.

\renewcommand\thefootnote{\dag}%

In the study of linear partial differential equations the use
of enlarged systems appeared in Grushin's work  \cite{Gr} on 
hypoelliptic operators. In a different context they were used in 
the thesis of the first author \cite{SjTh} and the $ \pm $ notation
comes from there -- see \S \ref{aens} for an explanation in the 
context of linear algebra.
As is seen there it is essential that the system is non-self-adjoint\footnote{That
distinguishes it from the KKT (for Karush-Kuhn-Tucker) systems popular in 
numerical studies -- see for instance \cite{KTT} -- which seem to be 
related to \S \ref{semp} below.}. For that historical, if somewhat personal 
reason, we refer to the problem 
\begin{gather}
\label{eq:gr}
\begin{gathered}
\left\{ \begin{array}{ll} P u + R_- u_- & = v \\
                         R_+ u & = v_+ \end{array} \right.\\
P \; : \; H_1 \; \rightarrow H_2 \,,  \ \ \ R_- \; : \; H_- \rightarrow H_2 \,, \ \ \ 
R_+ \; : \; H_1 \; \rightarrow \; H_+ \,, 
\end{gathered}
\end{gather}
as a {\em Grushin problem}. If it is invertible, 
we call it {\em well posed} and we write its inverse as follows
\begin{equation}
\label{eq:hi}
\left( \begin{array}{l} u \\ u_- 
\end{array} \right) = 
\left( \begin{array}{ll} E & E_+  \\ E_- & E_{-+} 
\end{array} \right)  
\left( \begin{array}{l} v \\ v_+ 
\end{array} \right) \,.
\end{equation}
In this case we will refer to $ E_{-+} $ as the {\em 
effective Hamiltonian} of $ P $.  That effective Hamiltonian normally has its
own physical interpretation as will be seen in examples in 
\S\S \ref{sebv}, \ref{aeps}, and \ref{aehf}. 

To illustrate this by a
straighforward example consider an operator $ P : L^2  ( \RR^n ) 
\rightarrow L^2 ( \RR^n ) $ defined as a convolution, $ P u = K \star u $, with
$ \widehat K \in L^\infty ( \RR ^n ) $. We can take 
$ H_\pm = L^2 ( \RR^n) $ and put $ R_-u_-( x )  = - (2 \pi)^{-n} \widehat u_-
( - x )  $, the negative of the inverse Fourier transform, and 
$ R_+ u (\xi ) = \widehat  u ( \xi ) $. One easily checks that the 
resulting Grushin problem is well posed and that $ E_{-+} $ is given 
by multiplication by $ \widehat K $. This of course is the effective 
Hamiltonian for the convolution operator which is invertible on $ L^2 $
if and only if $ \widehat K ^{-1} \in L^\infty $.

The main difficulty in constructing useful Grushin problems is the choice of
suitable operators $ R_\pm $ and of the spaces on which they act. As will 
be illustated below that depends on the situation even though one can notice
some underlying principles.

The paper is organized as follows. In \S \ref{se} we present in detail
several simple examples showing different ways of constructing Grushin problems.
In \S \ref{bt} we review basic linear algebra techniques which are useful when
studying Grushin problems arising in spectral theory. In \S \ref{sest}
we also show a typical parameter dependent estimate. Trace formul{\ae} which
are central in the study of classical/quantum correspondence are the 
subject of \S \ref{tf}: we give the basic idea in the context of Grushin 
problems and use it to prove the Poisson summation formula, in a way which
lends itself to many generalizations. Finally, in \S \ref{ae} we describe
-- without proofs -- 
four advanced examples: a remark on Lidskii-Lusternik-Vishik perturbation 
theory for matrices \cite{Lid}, \cite{MBO}, 
the Peierls substitution of solid state physics 
(from the work of Helffer and the first author \cite{HeSj1}), the quantum 
monodromy approach to the Gutzwiller trace formula, and  the
asymptotics of scattering poles in electromagnetic scattering by 
convex bodies (from earlier work of the authors \cite{SjZw},\cite{SjZw0}). It would
be very hard to survey all the examples in which the Grushin problem
appears explicitely -- not to mention, those in which it appears 
implicitely -- and we again made some personal choices.

\medskip
\noindent
{\sc Acknowledgments.} 
We would like to thank Steve Zelditch for suggesting a talk 
on Grushin problems
during the semi-classical semester at MSRI: this paper is a direct 
result of that. We are also grateful to Michael Overton for
the references to Lidskii's perturbation theory which are the basis
of \S \ref{lr}.

The work of the second author was supported in part by
the National Science Foundation  under the grant  DMS-0200732. 
He is also grateful to Universit\'e de Paris-Nord,
for its generous hospitality in October 2003.

\section{Simple examples}
\label{se}

We give five examples. The first two are  purely linear algebraic and
the third and fifth 
are intended to show how well known objects in mathematical physics 
fit in the Grushin problem
set up. The fourth example relates Grushin problems to analytic Fredholm 
theory which is one of the basic tools of spectral theory.

\subsection{The Moore-Penrose pseudoinverse}
\label{semp}
If $ P : \C^n \rightarrow \C^m $
is a linear transformation, its
{Moore-Penrose pseudoinverse} is the {\em unique} transformation
$ P^+ : \C^m \rightarrow \C^n $ satisfying
\begin{gather}
\label{eq:mp}
\begin{gathered}
P \circ P^+ \circ P = P \,, \ \ P^+ \circ P \circ P^+ = P^+\,,\\
( P \circ P^+)^* =  P \circ P^+
 \,, \ \ ( P^+ \circ P )^* = P^+ \circ P\,.
\end{gathered}
\end{gather}
If $ P $ has full rank then 
\[ P^+ =  \left\{ \begin{array}{ll} ( P^* P)^{-1} P^*  & n \leq m \\
P^* ( P P^* )^{-1} & n \geq m \end{array} \right.  \,.\]
In general $ P^+$ can be expressed by using the standard 
singular value decomposition
$ P = U^* \Sigma V $, and inverting the {\em non-zeros} entries in $ \Sigma$.
It is closely related to {\em least square problems} -- see 
\cite[Lecture 11]{BaTr}.

Another way to describe the pseudoinverse is as
\[ P^+ = \left( P \rest_{\ker ( P )^\perp } \right)^{-1} \pi_{\im ( P ) } \,, \]
where $ \pi_V $ is the orthogonal projection on the subspace $ V$, since
$ P\rest_{ \ker ( P ) ^\perp } : \ker ( P )^\perp \rightarrow \im ( P ) $ 
is bijective.

The pseudoinverse is a special case of $ E $ in \eqref{eq:hi}, with 
$ H_1 = \C^n $, $ H_2 = \C^m $, and for a natural choice of
$ R_\pm $, related to the least squares method. 
Before describing it, let us give
a general statement relating the Grushin problem to \eqref{eq:mp}:
\begin{prop}
\label{p:mp}
In the notation of \eqref{eq:gr} and \eqref{eq:hi} we always have
\[   E P E = E \,,\]
and the following equivalences
\[ \begin{split}
 P E P = P \ & \Longleftrightarrow E_- P = 0 \\
(PE)^* = PE \ & \Longleftrightarrow  ( R_- E_ -) ^* = R_- E_-  \\
(EP)^* = EP \ & \Longleftrightarrow  ( E_+ R_+ )^* = E_+ R_+  \;.
\end{split}\]
In particular, when the conditions on the left hold, $ E = P^+ $,
in the sense that equations in \eqref{eq:mp} are satisfied.
\end{prop}

We can now choose $ R_\pm $ so that the conditions in Proposition
\ref{p:mp} are satisfied. For that we simply put
\begin{gather*} 
H_- = P(H_1)^\perp \,, \ \ H_+ = \ker(P) \\
R_- : H_- = P(H_1)^\perp \hookrightarrow H_2 \,, \ \ R_+ : H_2 \stackrel{\perp}{\rightarrow}  
\ker(P) = H_+\,.
\end{gather*}

This is generalized in \S \ref{pse} in order to take into account 
small eigenvalues of $ (P^* P)^{\frac12} $ and  $ (P P^*)^{\frac12} $.

\subsection{Non-self-adjoint eigenvalue problems}
\label{aens}
In the previous example we had $ R^*_+ = R_- $. For non self-adjoint 
problems that constitutes an unnatural restriction as shown by the
following  elementary example. 

Let $ J $ be the $ n \times n $ upper triangular 
Jordan matrix:
\[ J = (J_{ij})_{1 \leq i,j \leq n } \,, \ \  J_{ij} = \left\{ 
\begin{array}{ll} 1 & \text{ for $k = j+1$ } \\
0 & \text{ otherwise. } \end{array} \right. \]

Let 
\[  e_+ = \left
( \begin{array}{l} 1 \\ 0 \\ \vdots \\ 0 \end{array} 
\right) \,, 
 \ \ e_- = \left( \begin{array}{l} 0 \\  \vdots \\ 0 \\ 1 \end{array} 
\right) \,.\]
Then $ J e_+ = 0 $, $ J^* e_- = 0 $, $ \| e_\pm \| = 1 $, and we can set 
up the following well posed Grushin problem for $ J - \lambda $:
\[   {\mathcal J} ( \lambda ) = \left( \begin{array}{ll} J - \lambda & R_- \\ 
R_+ & 0 \end{array} \right) \; : \;\
\CC^n \oplus \CC \longrightarrow 
\CC^n \oplus \CC \,, \ \ R_- u_- = u_- e_- \,, \ R_+ u = \langle u, e_+ \rangle.
\]
One easily checks that $ E_{-+} ( \lambda ) = \lambda^n $, and that 
\begin{gather}
\label{eq:epm}
\begin{gathered}  
E_+ v_+ = v_+  e_+( \lambda ) \,, \ \ 
E_- v = \langle v , e_- ( \lambda ) \rangle \,, \ \
  e_+ ( \lambda ) =  \left( \begin{array}{l} \ \ 1 \\ \ \ \lambda \\ \ \ 
\vdots \\ \lambda^{n-1}  \end{array} 
\right) \,, 
 \ \ e_- = \left( \begin{array}{l} \lambda^{n-1} \\ \ \ \vdots \\ \ \ \lambda \\ \ \ 1 \end{array} 
\right) \,.\end{gathered}
\end{gather}

\begin{figure}[htbp]
\begin{center}
\includegraphics[width=3.5in]{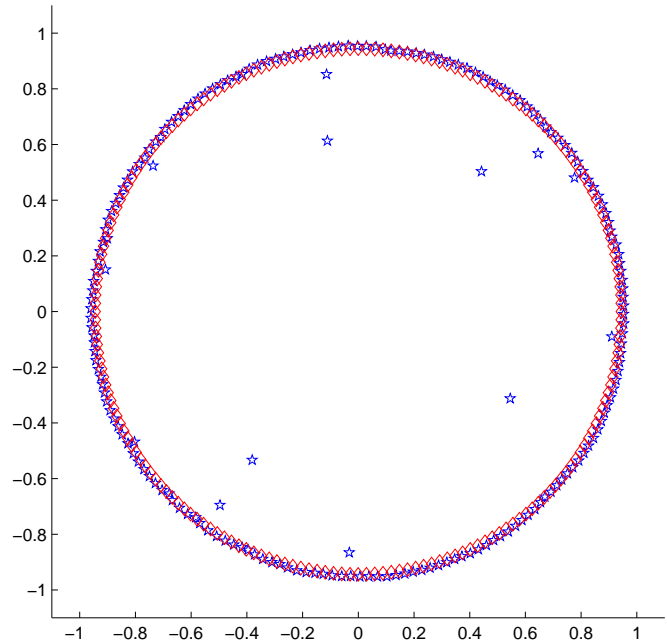}
\end{center}
\caption
{
\label{z}
Eigenvalues of a small random perturbation of a $ 200\times200$ Jordan block matrix (blue), 
and of the perturbation $ Q = \epsilon e_- \otimes e_+ $ (red)}
\end{figure}

If 
we add a small matrix perturbation, 
$ \epsilon Q $ to $ J$ the same problem remains well
posed and, using a Neumann series argument for matrices (see \S \ref{bfga} 
for a similar argument),
\[ E_{-+}^\epsilon ( \lambda ) = E_{-+} ( \lambda ) + \sum_{k=1}^\infty 
(-1)^k \epsilon^k E_- (\lambda ) Q ( E( \lambda ) Q )^{k-1} E_+ ( \lambda ) 
  \,,
\]
with uniform convergence for $ |\lambda | \leq \theta < 1 $ and  $
\epsilon \leq \epsilon_0 $, for some $ \epsilon_0 >0 $.
Using \eqref{eq:epm} we consequently see that
\begin{equation}
\label{eq:empe}  
E_{-+}^\epsilon ( \lambda ) = \lambda^n - \epsilon \langle Q e_+ 
( \lambda ) , e_- ( \lambda ) \rangle 
 + {\mathcal O} ( \epsilon^2   ) \,. \end{equation}
Hence when $ n $ is large and $ |\lambda | < 1 $ there will be no 
spectrum near $ \lambda $ for a generic perturbation $ Q $. This is
illustrated in Figure \ref{z}.
The most dramatic perturbation is 
obtained by taking $ Q $ with a large inner product 
$ \langle Q e_+ , e_- \rangle $.

This example is a linear algebraic  model of the first author's thesis 
\cite{SjTh} where the $ -+ $ notation was introduced. It was motivated by the
sign in H\"ormander's commutator condition -- see \cite{Sj} and also
\cite{ZwPse} for a lighthearted introduction. It is reflected here by the
fact that
\begin{equation}
\label{eq:-+}
[ J, J^*] e_\pm = \pm e_\pm \,. \end{equation}
This example will be revisited in a more general context in \S \ref{sest}.

\subsection{Feshbach method}

\label{aefr}
The Feshbach method which has been useful in the study of quantum resonances 
fits in the framework of Grushin problems discussed in this paper. 
To review it we follow \cite{DJ} and we refer to that paper for pointers
to the vast literature on the subject.

Suppose that a Hilbert space $ {\mathcal H} $ can be written as a direct 
sum $ {\mathcal H} = {\mathcal H}^{\rm{v}} \oplus {\mathcal H}^{\bar{\rm{v}}} $,
and that the quantum Hamiltonian under investigation decomposes
under this splitting as 
\[ H = \left( \begin{array}{ll} H^{\rm{v}\rm{v}} & H^{{\rm{v}}{\bar{\rm{v}}}} 
\\ H^{\bar{\rm{v}} \rm{v}} & H^{\bar{\rm{v}}\bar{\rm{v}}} \end{array} \right) \,.\]
Assume now that for $ z \in \Omega $, an open set in $ \CC$, the operator
$ ( z \bbbone^{\bar{\rm{v}} \bar{\rm{v}}} - 
H^{\bar{\rm{v}}\bar{\rm{v}}} ) $ is invertible.
Following \cite{DJ} we define the {\em resonance function}
\[ G_{\rm{v} } ( z ) = z \bbbone^{\rm{v}\rm{v}} - H^{\rm{v}\rm{v}} 
- H^{\rm{v}\bar {\rm{v}}} ( z \bbbone^{\bar{\rm{v}} \bar{\rm{v}}} - 
H^{\bar{\rm{v}}\bar{\rm{v}}} )^{-1} H^{\bar{\rm{v}}\rm{v}} \,, \]
which, in classical terminology reviewed 
in \S \ref{in} is just the {\em Schur complement} of
$  z \bbbone^{\bar{\rm{v}}\bar{\rm{v}}} - H^{\bar{\rm{v}}\bar{\rm{v}}} $ in $ z - H $.

It then follows, by block Gaussian elimination, 
that for $ z \notin \sigma (  H^{\bar{\rm{v}}\bar{\rm{v}}} ) $ 
\[ z \in \sigma ( H ) \ \Longleftrightarrow \ 0 \notin \sigma ( G_{\rm{v}} ( z ) ) \,,
\]
and moreover it can be verified directly that
\begin{equation}
\label{eq:fesh}
\tr \int_{\gamma_z } 
( \zeta - H ) ^{-1} d\zeta  = \tr  \int_{ \gamma_0}
\partial_\zeta G_{\rm{v}} ( \zeta ) G_{\rm{v}} ( \zeta )^{-1} d \zeta \,,\ \ 
\gamma_z (t ) = z +  \epsilon e^{it} , 0 \leq t \leq 2 \pi  \,,
\end{equation}
that is, the multiplicities agree.

To see how the Schur 
complement, and hence also the Feshbach method, fit in the Grushin
scheme we consider the following larger operator
\[ {\mathcal P}(z) = \left( \begin{array}{ll} z - H & R_- \\
\ \ R_+ & 0 \end{array} \right) \; : \; {\mathcal D} \oplus {\mathcal H}^{\rm{v}}
\longrightarrow  {\mathcal H} \oplus {\mathcal H}^{\rm{v}}  \,, \ \  
R_+ = ( \bbbone^{\rm{v}\rm{v}}
 \ {\bf 0}^{{\rm{v}}\bar{\rm{v}}} ) \,, \ 
R_- = \left( \begin{array}{l} \bbbone^{\rm{v}\rm{v}} \\
 {\bf 0}^{\bar{\rm{v}}{\rm{v}}} \end{array} \right)  \,. \]
If $ z \bbbone^{\bar{\rm{v}}\bar{\rm{v}}} -   H^{\bar{\rm{v}}\bar{\rm{v}}} $ is
invertible then this problem is well posed and Gaussian elimination shows that
\[ E_{-+} ( z ) =  - ( z \bbbone^{\rm{v}\rm{v}} - H^{\rm{v}\rm{v}} )
+  H^{\rm{v}\bar {\rm{v}}} ( z \bbbone^{\bar{\rm{v}} \bar{\rm{v}}} - 
H^{\bar{\rm{v}}\bar{\rm{v}}} )^{-1} H^{\bar{\rm{v}}\rm{v}} = 
- G_{\rm{v}} ( z ) \,.\]
The multiplicity formula  follows from general principles described in 
\S \ref{tbt} but of course it is easy enough to verify directly.

\subsection{Analytic Fredholm theory}
\label{seft}

Here we recall the discussion of the appendix in \cite{HeSj}. For the
basic facts from functional analysis we
refer to \cite{gk} for an in-depth treatment and to
\cite[Sect.19.1]{Horb} for a comprehensive introduction.

A bounded operator $ P : H_1 \rightarrow H_2 $ between two 
Banach spaces, is called a Fredholm operator if the kernel of $ P $,
\[  \ker P \defeq \{ u \in H_1 \; : \; P u = 0 \} \,, \]
and the cokernel of $ P$,
\[  \cker P \defeq H_2 / \{ Pu \in H_2 \; : \; u \in H_1 \} \,,\]
are finite dimensional. It then automatically follows (see for instance
\cite[Lemma 19.1.1]{Horb} or the comment after the proof of Proposition \ref{p:fred1}) 
that $ P H_1 $ is closed. For Fredholm operators
the index is defined as
\[ \ind P = \dm \ker P - \dm \cker P \,.\]
We have the following
\begin{prop}
\label{p:fred1}
Suppose that that for some choice of $ R_\pm $ the
Grushin problem \eqref{eq:gr} is well posed. Then $ P : H_1 \rightarrow H_2 $
is a Fredholm operator if and only if $ E_{-+} : H_+ \rightarrow 
H_- $ is a Fredholm operator, and 
\[ \ind P = \ind E_{-+} \,. \] 
\end{prop}
\begin{proof}
As for all well posed Grushin problems we have that $ R_+$, $ E_-$ are
surjective, and $ E_+ $, $ R_- $ are injective.

The equation $ P u = v $ is equivalent to 
\begin{equation}
\label{eq:sur} 
u = Ev + E_+ v_+ \,, \ \ 0 = E_- v + E_{-+} v_+ \,, 
\end{equation}
where $  v_+ = R_+ u $. This means that
\[ E_- \; : \; \im P \; \longrightarrow \im E_{-+} \,,\]
is surjective, and that it induces a bijective map
\[ E^\sharp_- \; : \; H_2/\im  P \; \longrightarrow H_- / \im  E_{-+} \,,\,.\]
In fact, if $ E_- v \in \im  E_{-+ } $  then we use  \eqref{eq:sur} to
see that $ v \in \im  P $. 

On the other hand,
\[ E_+ \; : \; \ker E_{-+} \; \longrightarrow \ker P \]
is a bijection. In fact, if $ u \in \ker P $ then $ u = E_+ v_+ $
and $ E_{-+} u_- = 0$ and the map is onto, which is all we need
to check as $ E_+ $ is always injective.

We conclude that 
\begin{equation}
\label{eq:dim}
 \dim \cker P = \dim \cker E_{-+} \,, \ \ \dim \ker P = \dim \ker E_{-+} 
\,.\end{equation}
In particular the indices are equal.
\end{proof}

For Fredholm operators we can always take $ H_\pm $ to be finite dimensional:
let $ n_+ = \dim \ker P $ and $ n_- = \dim \cker P $ and choose 
\[ R_- \; : \; \CC^{n_-} \; \longrightarrow H_2 \,, \ \ 
R_+ \; : \; H_1 \; \longrightarrow \; \CC^{n_+} \,,\]
of maximal rank and such that 
\[  R_- ( \CC^{n_-} ) \cap \im P = \{ 0 \} \,, \ \  \ker ( R_+ \rest_{\ker P} ) = \{ 0 \} 
\,.\]
In that case $ E_{-+} \; : \; \CC^{n_-} \rightarrow \; \CC^{n_+} $ and its index is,
of course, $ n_+ - n_- $. This argument also shows that the index does not change
under continuous Fredholm deformations of $ P $, and that $ P H_1 $ is closed:
by Banach's open mapping theorem the operators $ E_\bullet $ in \eqref{eq:hi}
(constructed using linear algebra only) are continuous.

The following standard result is proved particularly nicely using the
Grushin problem framework:
\begin{prop}
\label{p:4'}
Suppose that for $ z \in \Omega \subset \CC $, a connected open set,
$ A ( z  ) $ is a family 
of Fredholm operators depending holomorphically on $ z $. If $ A ( z_0 )^{-1} $
exists at a point $ z_0 \in  \Omega $. Then $ \Omega \ni z \mapsto A( z )^{-1} $is a meromorphic family of operators. 

\end{prop}
\begin{proof}
Let $ z_0 \in \Omega $ and let $ V ( z_0 ) $ be a small neighbourhood of $ z_0 $.
We can then form a Grushin 
problem for $ P = A ( z_0 ) $ as described before the statement of the proposition.
The same $ R^{z_0}_\pm $ 
give a well posed Grushin problem for $ P = A ( z ) $ for $ z\in V( z_0) $,
if $ V ( z_0 ) $ is sufficiently small. Since the index $ A ( z ) $ is equal to zero 
we see that $ n_+ = n_- = n $ and $ E^{z_0}_{-+} ( z ) $ is an $ n \times n $ matrix
with holomorphic coefficients. The invertibility of $ E^{z_0}_{-+} ( z ) $ is 
equivalent to the invertibility of $ A ( z ) $. 

This shows that there exists a locally finite covering of $ \Omega $,
$ \{ \Omega_j  \} $,
such that for $ z \in \Omega_j $, $ A ( z ) $ is invertible precisely
when $ f_j ( z ) \neq 0 $, where $ f_j $ is holomorphic in $ \Omega_j$.
Since $ \Omega $ is connected and since $ A ( z_1 ) $ is invertible for at
least one $ z_1 \in \Omega $ shows that {\em all} $ f_j$'s are not
identically zero. 

That means that $ \det E_{-+} ( z ) $ is
non-vanishing holomorphic function in $ V( z_0 ) $ and consequently $ E_{-+} ( z )^{-1} $
is a meromorphic family of matrices. Applying \eqref{eq:schur1} we conclude that
\[ A (  z)^{-1} = E( z ) - E_{+} ( z ) E_{ -+} ( z )^{-1} E_- ( z ) \]
is a meromorphic family of operators in $ V ( z_0 ) $, and since $ z_0 $ was arbitrary, in 
$ \Omega $.
\end{proof}

\subsection{Boundary value problems}
\label{sebv}

Let $ P $ be an elliptic second order operator on a compact manifold, 
$ X $,  with  
an orientable smooth boundary, $ \partial X $. For the simplest example
we could take $ P = - \partial_x^2 + V( x ) $ on $ [ a, b ] $, in which
case all the objects below are easily described.

We want to pose a Grushin problem for the Dirichlet realization of 
$ P$:
\[ P_D u = f \in L^2 ( X ) \,, \ \ u \rest_{\partial X } = 0 \,.\]

We then put
\[  \HH = H^2 ( X ) \cap H_0^1 ( X) \,, \ \ \HH_- = H^{\frac32} ( \partial  X) 
\,, \ \ \HH_+ = H^{\frac12} ( \partial X ) \,.\]
Let $ T \; : \; H^{\frac32} ( \partial X ) \rightarrow H^{2 } ( X ) $ be
an extension operator, with the following properties:
\[  T v \rest_{\partial X} = v  \,, \ \ \partial_\nu T v \rest_{\partial X } = 0 \,, \]
where $ \partial \nu $ is the outward normal differentiation at
$ \partial X $.
The operator $ T $ can, for instance, be obtained by introducing normal geodesic 
coordinates $ ( x , y ) $  in a collar neighbourhood of 
$ \partial X $, $ y \in \partial X $, and putting
\[  T v ( x , y ) = \chi ( x ) \exp (  x^2 \Delta_{\partial X} ) v (y) \,, \ \ 
- \Delta_{\partial X} \geq 0 \,, \]
where $ \chi \in \CIc ( [0, \delta) ) $, $ \chi \equiv 1 $, near $ 0 $. 

We then define
\begin{gather}
\label{eq:grub}
\begin{gathered}
  R_- \; : \;  \HH_- \longrightarrow \HH \,, \ \ 
R_+ \; : \; \HH \longrightarrow \HH_+ \\
R_- u_- \stackrel{\rm{def}}{=}  - P  T u_- \,, \ \ 
R_+ u \stackrel{\rm{def}}{=} \partial_\nu u \rest_{\partial X } \,.
\end{gathered}
\end{gather}

If we denote by $ P_N $ the Neumann realization of $ P $,
 \[ P_N u = f \in L^2 ( X ) \,, \ \ \partial_\nu u 
\rest_{\partial X } = 0 \,,\]
we have 
\begin{prop}
With $ R_\pm $ defined by \eqref{eq:grub} the Grushin problem for
$ P_D $ is well posed when $ P_N^{-1} $ exists. The effective 
Hamiltonian is given by the Neumann-to-Dirichlet map:
\[ E_{ - + } = N \,, \ \ N \; : \; \partial_\nu u \rest_{
\partial X }  \longmapsto u\rest_{\partial X} \,, \ \ P  u = 0 \,,\]
where the existence of $N$ is guaranteed by the invertibility of 
$ P_N $.
\end{prop}
\begin{proof}
We can write \eqref{eq:hi} explicitely using the Green operator,
$ G_N \stackrel{\rm{def}}{=} P_N^{-1} $, and the Poisson operator:
\[ Q_N f = u \,, \ \ P u = 0 \,, \ \ \partial_\nu u \rest_{\partial X } 
= f \,. \]
It can be easily constructed from $ P_N^{-1} $. 

Using this notation we have
\[\begin{split}
& E v  = G_N v  +  T ( (G_N v)\rest_{\partial X} )   \\ 
& E_- v =  (G_N v)\rest_{\partial X} \\
& E_+ v_+ =  Q_N v_+  +  T \circ N v_+  \\ 
&  E_{-+} v_+ = (Q_N v_+)\rest_{\partial X }  \end{split} 
\]
A direct verification proves the surjectivity. To 
prove injectivity we
see that injectivity of $ P_N $ gives 
\[ P ( u - T u_- ) = 0 \,, \ \ \partial_\nu ( u - T u_- ) = 0 \ 
\Longrightarrow \ u - T u_- = 0 \,. \]
Since $ u \rest_{\partial X } = 0 $  this
shows that $ u _- = T u_- \rest_{\partial X }  = 0 $, and hence $ u = 0$,
as well.
\end{proof}

A more systematic approach and one related to another use of two-by-two 
systems \cite{BdM},\cite[Sect.20.4]{Horb} can be described as follows. Suppose that 
$ P : \CI ( X ) \rightarrow \CI ( X ) $ is an elliptic operator of order $ m$,
and that we have two sets of boundary differential operators, with transversal 
orders $ < m $, 
\[ \begin{split}
& B_j \; : \; \CI ( X ) \rightarrow \CI ( \partial X ) \,, \ \  j=1,\cdots, J \,,\\
& C_k \; : \; \CI ( X ) \rightarrow \CI ( \partial X ) \,, \ \  k=1,\cdots, K \,. 
\end{split}
\]
For instance we can consider $ P = \Delta $, 
$ B_1 u =  \partial_\nu u \rest_{\partial X} $,
$ C_1 u = u \rest_{\partial X } $, $ J=K=1$.

We want to study the boundary problem 
\begin{equation}
\label{eq:obl}
P u = f \ \ \text{in $ X $}, \ \ C_k u = h_k \ \ \text{in $ \partial X $}, \ \ 
k=1,\cdots, K \,, \end{equation}
assuming that
the boundary problem 
\begin{equation}
\label{eq:wp}
 P u = f \ \ \text{in $ X $}, \ \ B_j u = g_k \ \ \text{in $ \partial X $}, \ \ 
j=1,\cdots, J \,, \end{equation}
is well posed.
To avoid technical issues involving Sobolev spaces (see \cite[Chapter 20]{Horb}) we
will remain in the $ \CI $ category.
We then put:
\[ {\mathcal H}_1 = \CI ( X) 
\,, \ \ {\mathcal H}_2 = \CI ( X) \otimes \CI ( \partial X )^K \,, 
\ \ 
{\mathcal H}_{-}  = \CI ( \partial X )^K \,, \ \  {\mathcal H}_{+} = \CI ( \partial X)^J  \,, \]
writing
\[ u_- = \left( \begin{array}{l}  u_-^1  \\ 
\ \vdots \\ u_-^K \end{array}  \right) \in   {\mathcal H}_-  \, , \ \ 
 v = \left( \begin{array}{l}  v_X  \\ v_{\partial X }^1 \\ 
\ \vdots \\ v_{\partial X}^K \end{array}  \right) \in   {\mathcal H}_2 \,, \]
and define
\begin{gather}
\label{eq:elbvp} 
\begin{gathered}
Q u \defeq \left( \begin{array}{l}
P u \\ C_1 u \\ \ \ \vdots \\ C_K u \end{array} \right) \,,
 \ \ 
  R_- u_- \defeq \left( \begin{array}{l} \ 0 \\ u_-^1 \\ \ \vdots \\ u_-^
K  \end{array} \right) \,, \ \ 
R_+ u  \defeq \left( \begin{array}{l} B_1 u \\ \ \ 
\vdots \\ B_J u   \end{array} \right) \,. \end{gathered} \end{gather}
We have the following {\em formal} 
\begin{prop}
\label{p:bvp}
Suppose that the boundary value problem \eqref{eq:wp} is well posed. Then the 
Grushin problem 
\[   Q u + R_- u_- = v \,, \ \ R_+ u = v_+ \,,\]
obtained using the operators \eqref{eq:elbvp} is well posed and the effective
Hamiltonian, 
\[ E_{-+} \; : \; \CI(\partial X)^J \; \rightarrow \; \CI( \partial X)^K \]
is a generalization of the Neumann-to-Dirichlet map:
\begin{equation}
\label{eq:fredb}
 E_{-+}  \; : \; \left( \begin{array}{l}  v_+^1 \\ \ \vdots \\ v_+^
J  \end{array} \right)  \longmapsto 
\left( \begin{array}{l} 
 C_1 u \\ \ \ \vdots \\ C_K u \end{array} \right) \,, \ \
P u = 0 \,, \ \  B_j u = v_+^j \,, \ \ j = 1, \cdots, J \,.\end{equation}
\end{prop}

\medskip

For boundary value problems one of the basic issues is showing that, on suitably 
chosen spaces, the operator $ u \mapsto (Pu, C_1u , \cdots C_k u )$ has
the Fredholm property. By Proposition \ref{p:fred1} that is equivalent to showing the
Fredholm property of the operator \eqref{eq:fredb}. The reduction to the boundary 
described in Proposition \ref{p:bvp} will furnish us with another example in 
\S \ref{bfip}.

\section{Basic techniques}
\label{bt}

Here we present some general results about systems arising from 
considering {\em Grushin problems} and examples showing how they can be used.
 We recall that a Grushin problem
for an operator $ P : H_1 \rightarrow H_2 $ is a system
\begin{equation}
\label{eq:gr'}
\left\{ \begin{array}{ll} P u + R_- u_- & = v \\
                         R_+ u & = v_+ \end{array} \right.
\end{equation}
where $ R_- : H_- \rightarrow H $, $ R_+ : H \rightarrow H_+ \,, $.
In matrix form we can write
\[ {\mathcal P} \defeq \left( \begin{array}{ll} \ P & R_- \\
R_+ & \ 0 \end{array} \right) \; : \; H_1 \oplus H_- \; \longrightarrow 
H_2 \oplus H_+ \,. \]
We say that the Grushin problem is {\em well posed} if we have 
the inverse
\[ {\mathcal E} = \left( \begin{array}{ll} E & E_+  \\ E_- & E_{-+} 
\end{array} \right)   
  \; : \; H_2 \oplus H_+ \; \longrightarrow 
H_1 \oplus H_- \, ,\]
that is
\begin{equation}
\label{eq:hi'}
\left( \begin{array}{l} u \\ u_- 
\end{array} \right) = 
\left( \begin{array}{ll} E & E_+  \\ E_- & E_{-+} 
\end{array} \right)  
\left( \begin{array}{l} v \\ v_+ 
\end{array} \right) \,.
\end{equation}
In this case we will refer to $ E_{-+} $ as the {\em 
effective Hamiltonian} of $ P $.

\subsection{Two by two systems}
\label{tbt}

Here we consider an invertible system
\begin{gather}
\label{eq:tts}
\begin{gathered}
{\mathcal A} \defeq \left( \begin{array}{ll} A_{11} & A_{12} \\
A_{21} & A_{22} \end{array} \right) \; : \; H_{1} \oplus H_{2} \; \longrightarrow 
\; \tH_1 \oplus \tH_2 \,, \\ 
{\mathcal B} \defeq {\mathcal A}^{-1} = 
\left( \begin{array}{ll} B_{11} & B_{12} \\
B_{21} & B_{22} \end{array} \right) \; : \; \tH_{1} \oplus \tH_{2} \; \longrightarrow 
\; H_1 \oplus H_2 \,.
\end{gathered}
\end{gather}

We first recall the formula involving an expression 
known as the  {\em Schur complement} in linear 
algebra and as the {\em Feshbach operator} in mathematical physics:
\begin{lem}
\label{l:schur}
Suppose that $ A_{22} $ is invertible. Then $ B_{11} $ is invertible, and 
\begin{equation}
\label{eq:schur}
B_{11}^{-1} =  A_{11} - A_{12} A_{22}^{-1} A_{21} \,.
\end{equation}
\end{lem}
\begin{proof}
Using $ B_{11} A_{11} + B_{12} A_{21} = I $ and $ B_{11} A_{12} = -B_{12} A_{22}$
we see that
\[ B_{11} A_{11} - B_{11} A_{12} A_{22} ^{-1} A_{21} = 
I - ( B_{12} + B_{11} A_{12} A_{22}^{-1} ) A_{21} = I \,,\]
and the left inverse property is derived similarly.
\end{proof}

We now allow the entries of $ {\mathcal A} $ to depend on a parameter, and 
denote differentiation with respect to that parameter by $ A \mapsto \dot A $.
The next lemma explicitely shows that the traces of $ \dot B_{11} B_{11}^{-1} $
and $ \dot A_{22} A_{22}^{-1} $ differ by terms not involving any inverses. 
In the case of holomorphic dependence on the parameter that means that 
these traces differ by holomorphic terms which disappear in contour integration.
Before stating this precisely let us recall some basic facts about
trace class operators -- see \cite{gk} or \cite[Sect.19.1]{Horb}.

If $ H_j $ are infinite dimensional Hilbert spaces, 
the operator $ A : H_1 \rightarrow H_2 $
is said to be {\em of trace class} if the self-adjoint operator $ 
( A A^* )^{\frac12} : H_2 \rightarrow H_2 $ has a discrete spectrum,
 $ \{\mu_j\}_{j=1}^\infty  $, and $ \sum_{j=1}^\infty \mu_j < \infty $.
If $ A $ is of trace class, and $ B_1 \rightarrow H_1 \rightarrow H_1 $,
$ B_2 : H_2 \rightarrow H_2 $ are bounded operators then $ A B_1 $ and
$ B_2 A $ are of trace class. 

If $ H = H_1 = H_2 $, $ A $ is of trace class we can define the
{\em trace} of $ A $ as follows. Let 
$ \{e_j\}_{j=1}^\infty  $  be an orthonormal basis of $ H $, then
\[ \tr A \defeq \sum_{j=1}^\infty \langle A e_j , e_j \rangle_{H} \,,\]
and this definition is independent of the choice of a basis.

Finally, if $ A : H \rightarrow H $ is of trace class and $ B : H \rightarrow H $
is bounded then
\begin{equation}
\label{eq:ct}
\tr [ A , B ] = 0 \,, \ \ [ A , B ] = AB - BA \,.
\end{equation}

\begin{lem}
\label{l:2}
Suppose that the operators 
$ \dot A_{ij} $  are of trace class. 
Then, when $ A_{22} $ is invertible, we have
\begin{gather}
\label{eq:l2}
\begin{gathered}
 \tr B_{11}^{-1} \dot B_{11} = 
\tr A_{22}^{-1} \dot  A_{22} - \tr \dot {\mathcal A}  {\mathcal B } \,.
\end{gathered}
\end{gather}
\end{lem}
\begin{proof} 
This is a straighforward computation based on the formul{\ae}, $ 
{\mathcal A} {\mathcal B} = I$,  $ \dot {\mathcal B}
= - {\mathcal B} \dot {\mathcal A} {\mathcal B} $, cyclicity of the
trace, and Lemma \ref{l:schur} (we note that all $ \dot B_{ij} $, and in particular
$ \dot B_{11} $, are of trace class). We obtain 
\begin{gather*}
 \tr B_{11}^{-1} \dot B_{11} = 
\tr A_{22}^{-1} \dot  A_{22} + 
\tr E_1 + \tr E_2 \,, \\  
E_1 = - \dot A_{11} B_{11} - \dot A_{12} B_{21} \; : \; \tH_1 \longrightarrow \tH_1  \,, \\
E_2 =  -  \dot A_{21} B_{12} - \dot A_{22} B_{22}  \; : \; \tH_2 \longrightarrow \tH_2\,,  
\end{gather*}
and
\[  \tr  E_1 + \tr E_2 = - \tr \dot {\mathcal A}  {\mathcal B } 
= \tr {\mathcal A}  \dot {\mathcal B } \,. \]
\end{proof} 

The relevance of this discussion for Grushin problems (which in 
principle have $ B_{22} = 0 $) will become apparent in the next subsection.

\subsection{From one Grushin problem to another}
\label{fogto}

Suppose that we have a well posed Grushin problem \eqref{eq:gr} with the 
inverse given by \eqref{eq:hi}.

We want to check if another Grushin problem is well posed:
 \begin{equation}
\label{eq:tgr}
\left\{ \begin{array}{ll} P u +  \tR_-  \tu_- & = \tv \\
                         \tR_+ \tu & = \tv_+ \end{array} \right.
\end{equation}
The corresponding operator will be denoted by 
$ \widetilde {\mathcal P} \; : \; 
H_1 \oplus \tH_- \rightarrow H_2 \oplus \tH_+ $. 
If the inverse exists we will denote it by $ \widetilde {\mathcal E} $, with the
corresponding notation for the entries.

The simple answer is given in 
\begin{prop}
\label{p:gr}
The Grushin problem \eqref{eq:tgr} is well posed if and only if 
the following system of operators obtained from the 
solution \eqref{eq:hi} of the well posed problem \eqref{eq:gr},
\begin{equation}
\label{eq:redgr}
{\mathcal G} = \left( \begin{array}{ll}   - \tR_+ E 
\tR_-  &  \tR_+ E_+ \\   - E_- \tR_-   & \ \ E_{-+}  
\end{array} \right)  \,, \end{equation}
is invertible, that is if and only if the matrix of operators has a two sided 
inverse. In that case
\begin{equation} 
\label{eq:redtgr}
{\mathcal G}^{-1} = \left( \begin{array}{ll} 0 & I \\ I & 0 \end{array} 
\right) \widetilde {\mathcal G}  \left( \begin{array}{ll} 0 & I \\ I & 0 \end{array} 
\right) \,, \end{equation}
where 
\[ \widetilde 
{\mathcal G} = \left( \begin{array}{ll}   - R_+ \widetilde E  R_- &  R_+ \widetilde E_+ 
 \\ 
  - \widetilde E_- R_-    & \ \  \widetilde E_{-+}  
\end{array} \right)  \,.\]
\end{prop}
\begin{proof}
In place of \eqref{eq:tgr} we can consider a larger system 
\[\left\{ \begin{array}{ll} P u  + R_- u _- +   \tR_-  \tu_- & = \tv  \\
                         R_+ u & = v_+ \\
                         \tR_+ u & = \tv_+  \end{array} \right. \,, 
\]
in which $ \tv $, $ \tv_+ $, and $ u_- $ are given, and $ u $, $ \tu_- $, and $ v_+ $ 
are unknown. We can solve \eqref{eq:tgr} by putting $ u_- = 0 $.
Using 
\eqref{eq:hi} with we can write 
\begin{equation}
\label{eq:v+}
 \begin{split}
u & = E ( \tv - \tR_- \tu_- ) + E_+ v_+ \\
u_- & = E_- ( \tv - \tR_- \tu_- ) + E_{-+} v_+ \,, 
\end{split}\end{equation}
or, since $ \tR_+ u = \tv_+  $,
\[ \begin{split}
& \tR_+ E_+ v_+ - \tR_+ E \tR_- \tu_- = \tv_+ - \tR_+ E \tv \\
& E_{-+} v_+ - E_- \tR_- \tu_- = u_- - E_- \tv  \,.
\end{split}\]
which in turn 
can be rewritten as 
\[ {\mathcal G } \left( \begin{array}{l}  \tu_- \\ v_+ \end{array}
\right) = \left( \begin{array}{l} \tv_+ - \tR_+ E \tv \\ u_-  -  E_- \tv 
\end{array} \right) \,.\]
Hence the invertibility of $ {\mathcal G } $ implies that \eqref{eq:tgr} is well 
posed. In fact, we first obtain $ \tu_- $ by inverting $ {\mathcal G} $ and then 
$ u $ by using the first equation in \eqref{eq:v+}. When $ \tv = 0 $ we see that
\[ {\mathcal G } \left( \begin{array}{l}  \tu_- \\ v_+ \end{array}
\right) = \left( \begin{array}{l} \tv_+  \\ u_-   
\end{array} \right) \,,\]
from which the equivalence and \eqref{eq:redtgr} follow.
\end{proof}

We  illustrate Proposition \ref{p:gr}
with an example which is also the basis for \S \ref{tfpf}
below. Let us consider 
\begin{equation}
\label{eq:haci}
   P = P ( z ) \defeq hD_x - z \,, \ \ x \in \SP^1 \defeq \RR/( 2 \pi \ZZ ) \,.
\end{equation}
We formulate a Grushin problem as in \cite[Sect.2]{SjZw} where it was
motivated by \cite{HeSj2}. For that we want to find $ R_\pm ( z ) $ so that
\begin{equation}
\label{eq:2.gr}
{\mathcal P } ( z ) \stackrel{\rm{def}}{=} 
\left( \begin{array}{ll} P - z & R_- ( z) \\
R_+ ( z )  & 0 \end{array} \right) \; : \; H^1 ( {\mathbb S} ^1 ) \times \CC 
\ \longrightarrow \ L^2 ( {\mathbb S} ^1 ) \times \CC \,,
\end{equation}
is invertible. Rather than give the answer in a ``deus ex machina'' manner
we follow our original reasoning. 
First, a boundary condition 
\[  R_+ u \stackrel{\rm{def}}{=} u ( 0 ) \,, \]
is a natural choice. 
Then we can locally solve
\[  \left\{ \begin{array}{l} ( P - z ) u = 0 \\
R_+ u = v \end{array} \right.  \,,  \] 
by putting 
\[ u = I_+( z)  v  = \exp ( i z x / h ) v \,, \ \ \  - \epsilon < x < 
2\pi - 2 \epsilon \,.\]
This is the forward solution, and we can also define the backward one
by 
\[ u = I_- ( z) v  = \exp ( i z x / h ) v \,, \ \ \  - 2 \pi  + 2 \epsilon < x 
 \epsilon \,.\]
The monodromy operator $ M ( z , h ) : \CC \rightarrow \CC $, can be 
defined by 
\begin{equation}
\label{eq:2.mon}  I_+ ( z) v ( \pi )  = I_- ( z ) M ( z , h )v  ( \pi ) 
 \,, \end{equation}
and we immediately see that
\[ M ( z , h ) = \exp ( 2 \pi i  z /{ h }) \,. \]
We use $ I_\pm ( z) $ and the point $ \pi $ to work with objects 
defined on $ {\mathbb S}^1 $ rather than on its cover: a more intuitive
definition of $ M ( z , h ) $ can be given by looking at a value of the
solution after going around the circle.

Let $ \chi \in \CI ( {\mathbb S}^1 , [0 , 1 ] ) $  have the properties
\[ \chi ( x ) \equiv 1 \,, \ \ - \epsilon < x < \pi + \epsilon \,, \ \ 
\chi ( x) \equiv 0 \,, \ \   - \pi + 2 \epsilon < x < -2 \epsilon \,, \]
and put 
\[ E_+ ( z ) =  \chi I_+ ( z) + ( 1 - \chi ) I _- ( z ) \,.\]
We see that 
\[  ( P - z) E_+  =  [ P , \chi ] I_+ ( z) - [ P , \chi ] I_- ( z) 
=   [ P , \chi ]_- I_+ ( z) - [ P , \chi ]_- I_- ( z)  \,, \]
where $ [ P , \chi]_- $ denotes the part of the commutator supported near
$ \pi $.  This can be simplified using \eqref{eq:2.mon}:
\[  (i/h) ( P - z ) E_+ + ( i / h ) 
[ P , \chi] _- I_- ( z ) ( I - M ( z , h ) ) = 0 \,, \]
which suggests putting 
\[ R_- ( z ) = ( i / h ) [ P , \chi]_- I_- ( z) \,, \]
so that the problem
\[ \left\{ \begin{array}{l} ( P - z) u + R_- ( z) u_ - = 0 \\
R_+ ( z ) u = v \end{array} \right. \]
has a solution:
\[ \left\{ \begin{array}{l} u = E_+ ( z ) v \\
u_- = E_{- + } ( z) v  \end{array} \right. \,, \]
with $ E_{-+} ( z) = I - M( z , h )\,. $ In fact, 
it is much more natural, and easier for full-blown microlocal generalizations,
to consider a different $ R_+ ( z ) $
so that, with symmetry reminiscent of \S \ref{aens},
\begin{equation}
\label{eq:remi}
  R_- ( z ) u_- = u_-  e_- ( z ) \,, \ \ R_+ ( z ) u = \langle u , e_+ ( z ) 
\rangle\,, \ \  e_\pm ( z,x ) =  ( i / h ) [ P , \chi ]_{\pm }
( \exp( {i \bullet z}/{h} ) )
( x) \,.\end{equation}

One can show 
that with this choice of $ R_\pm ( z ) $, \eqref{eq:2.gr}
is invertible and then
\[ {\mathcal P}(z)^{-1} = {\mathcal E} ( z) = 
\left( \begin{array}{ll} E ( z ) & E_+ ( z) \\
E _- ( z ) & E_{-+ } ( z ) \end{array} \right) \,, \]
where all the entries are holomorphic in $ z $, and 
$ E_+ ( z) $, $ E_{ -+}( z) $, are as above. The operator $ E_{-+} ( z) $
is the effective Hamiltonian in the sense that its invertibility controls
the existence of the resolvent:
\begin{equation}
\label{eq:res}
 ( P - z )^{-1} = E ( z) - E_+ ( z) E_{ - + } ( z) ^{-1} E_- ( z) \,. 
\end{equation}

The invertibility is independent of $ \chi $ with the properties
described above. Hence we can  move to a singular limit in the
choice of $ \chi$ and deform $ \pi $ to $ 0 $. That means that 
we consider the following
Grushin problem (with suitably modified spaces):
\[ \left\{ \begin{array}{ll}  (i/h)( hD_x - z )  u ( x ) 
 - \delta_0 ( x )  u_- & = v ( x ) \\
                          u ( 0+) & = v_+ \end{array} \right. \]
In fact, we can write 
\[  u = E ( z ) v + E_+ ( z ) v_+ \,, \ \ u_- = E_- ( z ) v + E_-{ -+ } ( z ) v_+ \,,
\]
where 
\begin{gather}
\label{eq:ee}
\begin{gathered}
 E( z ) v ( x ) 
= \bbbone_{ [ 0 , 2 \pi [ }
\int_0^x \exp ( i ( x - y ) z/h ) v ( y ) dy \,, \ \ 
E_+ ( z ) v_+ = v_+ \exp ( i x z / h ) \bbbone_{[0, 2 \pi [ } \,, \\
E_- ( z ) v  =  - \exp ( 2 \pi i z /h ) \int_0^{2 \pi } 
\exp (  - i yz/h ) v ( y ) dy\,,  \ \
E_{-+} ( z ) = 1 -  \exp ( 2 \pi i z / h ) \,. \end{gathered}
\end{gather}

We finally come to an application of Proposition \ref{p:gr}.
In \eqref{eq:remi} it would be nice to be able to take $ e_- ( z ) = e_+ ( z )$,
that is to have a self-adjoint Grushin problem. That would also
simplify matters in more complicated situations. Hence suppose that
\[   e_\pm ( z , x ) = f ( x) \exp ( i x z / h ) \,. \]
We then have to consider the invertibility of the matrix $ {\mathcal G} $ 
in Proposition \ref{p:gr} -- which here is an honest $ 2 \times 2 $ matrix.
A brief calculation shows that for $ z \in \RR $, 
$ {\mathcal G} $  is equal to 
\[ \left( \begin{array}{ll} \ \  - A & \ \ \ \overline B \\
 e^{ 2\pi i z/ h } B & 1 - e^{2\pi i z / h } \end{array} \right) \,, \ \
A = \int_0^{2 \pi } \int_0^x  \overline {f ( x )} f ( y ) dy dx \,, \ \ 
B = \int_0^{2 \pi } f ( x ) dx \,,\]
and we observe that $ |B|^2 = A + \overline A $.

Hence the condition for invertibility becomes
\[ \Re ( A  e^{ - \pi i z / h }  )   \neq 0 \,, \]
and that will always be violated for some $ z \in \RR $. Hence we cannot
have a well posed Grushin problem for all $ z \in \RR $ with $ e_- = e_+ $ 
in \eqref{eq:remi}.

\subsection{Iterated problems}
\label{bfip}

The Grushin problems can be iterated and this is particularly 
important when the intermediate Grushin problems are {\em formal}
and only after one or more iterations we obtain a well posed
problem. An example of a useful formal problem will be given in 
\S \ref{aegt}.

Before giving an example of that we start with the following 
simple
\begin{prop}
Suppose that \eqref{eq:gr} is well posed with the inverse given by 
\eqref{eq:hi}. If 
\[  \left( \begin{array}{ll} E_{-+} & N_- \\
N_+ & \ 0 \end{array} \right) \; : \; H_+ \oplus V_- \; \longrightarrow
H_- \oplus V_+ 
\]
is invertible, with the inverse 
\[ {\mathcal F} = \left( \begin{array}{ll} 
F & F_+ \\ F_- & F_{-+} \end{array} \right) \,,\]
then the  new Grushin problem 
\[ \left( \begin{array}{ll} \ P & R_- N_- \\ N_+ R_+ & \ 0 \end{array} \right)
 \; : \; H \oplus V_- \; \longrightarrow
H \oplus V_+ 
\,,\] 
is well posed with the inverse given by 
\[ \left( \begin{array}{ll} E - E_+F E_- & E_+ F_+ \\
\ \ F_- E_- & \; -  F_{-+} \end{array} \right) \,.\]
\end{prop}
\begin{proof} 
We need to solve
\[ \begin{split} 
P u + R_- N_- u _- & = v \\
N_+ R_+ u & = v_+ \end{split} \]
Putting $ N_- u_- = \tilde u_- $, and $ R_+ u = \tilde v_+  $, we
obtain 
\[ \begin{split} 
P u + R_- \tilde u _- & = v \\
R_+ u & = \tilde v_+ \end{split} \]
which is solved by taking 
\[  \begin{split} 
&  u   =  E v  + E_+ \tilde v_+ \\
&  \tilde u_-  = E_- v + E_{-+} \tilde v_+ \end{split} \]
Recalling the definitions of $ \tilde u_- $ and $ \tilde v_+ $
this becomes 
\[ \left( \begin{array}{ll} E_{-+} & N_- \\
\ N_+ & \ 0 \end{array} \right) \left( \begin{array}{l} \; \  \tilde v_+ \\ - u_- \end{array}
\right) = \left( \begin{array}{l} - E_- v \\ \ \ \ 
v_+ \end{array} \right) \,.\]
Solving this using $ {\mathcal F} $ gives the lemma.
\end{proof}

We will mention one concrete example for which iterated Grushin problems
are useful. In the notation of Proposition \ref{p:bvp} consider for $ X $ 
an open set in $ \RR^{n+1} $, with a smooth boundary $ \Omega $, and
put $ P u  = \Delta u $, $ B_1 u = u\rest_{\Omega } $, and $ C_1 u 
= V u \rest_{\Omega } $, $ K=J=1$, where $ V $ is a vectorfield. 
If $ V $ is {\em not} everywhere transversal to $ \Omega $ we
obtain the {\em oblique derivative problem} and the operator 
\eqref{eq:fredb} is not elliptic and not self-adjoint. As in 
\cite{SjTh} one can then construct a new Grushin problem for that
operator using the structure of the set where $ V $ is not transversal to 
$\Omega $. A ``baby'' version of that type of problem
was presented on the level of matrices in 
\S \ref{aens}.

\subsection{A Grushin approximation scheme}
\label{bfga}

Let $ H $ be a Hilbert space, and $ H_0 $ a finite dimensional subspace
with an orthonormal basis $ \{ e_j\}_{j=1}^N $.
Let us introduce 
\[  R_+ \; : \; H \rightarrow \CC^N \,, \ \ 
R_- = R_+^* \; : \; \CC^N \rightarrow H \,, \]
given by 
\[ (R_+ u)_j  = \langle u , e_j \rangle \,, \ \ 
R_- u_- = \sum_{j=1}^N u_{-,j} e_j \,.\]
We want to consider the Grushin problem for the operator
\[ P = I - T \,, \ \ T \; : \; H \rightarrow H \,.\]
In many interesting situations we can reduce the study of 
a differential operator to the study of $ I - T $ by 
factoring out an invertible term.

The following lemma is related to the example presented in  \S \ref{semp}:
\begin{lem}
\label{l:4}
Let $ \pi $ be the orthogonal projection on the span of $ e_j$'s.
With the operators $ R_\pm $ given above, and $ P = ( 1 - \pi T ) $ the 
problem \eqref{eq:gr} is well posed, and the matrix \eqref{eq:hi}
is given by 
\begin{equation}
\label{eq:matrix}
\left( \begin{array}{ll} E^0 & E_+^0  \\ E_-^0 & E_{-+}^0 
\end{array} \right)  
=
\left( \begin{array}{ll} 
\ \ \ \ 1 - \pi  &  \ \ \ \ \ R_-   \\  R_+ ( I + T ( 1 - \pi ) )  &   
R_+ T R_- - 1 
\end{array} \right)  \,.
\end{equation}
\end{lem}
\begin{proof}
We observe that 
\begin{equation}
\label{eq:obs}
 R_- R_+ = \pi \,, \ \ R_+ R_- = \Id_{\CC^N} \,, \ \ \pi R_- = R_- \,, \ \ 
R_+ \pi = R_+ \,,\end{equation}
which leads to an immediate verification of \eqref{eq:matrix}.
\end{proof}

We can now consider the problem for $ 1 - T $ and we have
\begin{prop}
\label{p:3}
If $ \| ( 1  - \pi ) T  \| < \delta < 1 $ then 
the Grushin problem \eqref{eq:gr} with $ P = 1 - T  $ and $ R_\pm $ 
as in Lemma \ref{l:4}, is well posed, and the effective Hamiltonian 
has the following expansion:
\[ E_{- + } = E_{ -+}^0 + \sum_{ k=1}^\infty R_+ T ( ( 1 - \pi )T)^k 
R_- \,. \]
\end{prop}
\begin{proof}
This is a typical Neumann series argument. Using Lemma \ref{l:4}, 
and writing $ I - T = I - \pi T - ( I - \pi ) T $ we see that
\[
\begin{split}
&  \ \ \  \left( \begin{array}{ll} I - T & R_- \\
\ \ R_+ & 0 \end{array} \right) = 
\left( \begin{array}{ll} I - \pi T &  R_- \\
\ \ R_+ & 0   \end{array} \right)  \left( 
\Id_{ H \oplus \CC^N } - 
\left( \begin{array}{ll}  \ \ (I - \pi)T  & 0  \\
 R_+ T ( I - \pi ) T   &   0  \end{array} \right)  \right)  \,, \end{split} \]
where we used \eqref{eq:obs} to multiply
\[ \left( \begin{array}{ll} 
\ \ \ \ 1 - \pi  &  \ \ \ \ \ R_-   \\  R_+ ( I + T ( 1 - \pi ) )  &   
R_+ T R_- - 1 
\end{array} \right) \left( \begin{array}{ll} ( I - \pi ) T  & 0 \\
\ \ \ 0 & 0 \end{array} \right) \,.\]
Hence 
\[ \left( \begin{array}{ll} E & E_+  \\ E_- & E_{-+} 
\end{array} \right)   = \sum_{ k=0}^\infty \left( \begin{array}{ll}  \ \ (I - \pi)T  & 0  \\
 R_+ T ( I - \pi ) T   &   0  \end{array} \right)^k 
 \left( \begin{array}{ll} 
\ \ \ \ 1 - \pi  &  \ \ \ \ \ R_-   \\  R_+ ( I + T ( 1 - \pi ) )  &   
R_+ T R_- - 1 
\end{array} \right) \,.\]
Since 
\[ \left( \begin{array}{ll} A & 0 \\ B & 0 \end{array} \right)^k =
\left( \begin{array}{ll} \ \ A^k & 0 \\ B A^{k-1}  & 0 \end{array} \right) \,,
\]
we immediately obtain the formula for $ E_{-+} $.
 \end{proof}

The difficulty with the approximation scheme described here is
the need for an orthonormal basis. In practice that is rarely given 
in theoretical and, especially, numerical problems. To some extend that
can be remedied as follows.


Let 
$ T $ 
be
an operator with the property replacing 
the smallness of $ ( I - \pi ) T $: 
we postulate that there exists
a finite set, $ \{ e_j\}_{j=1}^M $, 
with the following property: 
\begin{equation}
\label{eq:t} 
  \forall \; u \in H  \ \exists \; t_j \in \CC \,, r \in H \ \  T u = \sum_1^M t_j e_j + r \,, \ \ 
\| r \| \leq \delta \|u\| \,, \ \  \frac{1}{C} 
\| T u \| \leq \| \, \vec t \, \|_{\ell^2} \leq C \| T u \|\,. 
\end{equation}

As before we would like to construct a well posed (in the sense that
its stability constant is controlled, not just that it is 
invertible) Grushin problem for $ I - T $. 

First, we need to modify the spanning set, $ \{ e_j\}_{j=1}^M $.
For that we introduce the Grammian matrix,
\[ G \stackrel{\rm{def}}{=} ( \langle e_i, e_j \rangle)_{ 1 \leq i,j \leq M }
\,.\]
It is positive semi-definite and hence can be diagonalized. 
We then can, after a unitary (in $ \CC^M $) ``reorganization'', 
assume that $ \{ e_j \} $ satisfy 
\[  \langle e_i , e_j  \rangle=\delta_{ij} \lambda_j  \,, \ \ 
\lambda_1 \geq \lambda_2 \geq \cdots \geq \lambda_M \geq 0 \,.\]

Suppose that 
 $ \lambda_j > (\epsilon/C )^2 $ for $ j \leq L $. 
The condition number, $ \| G \| \| G^{-1} \| $, of the Grammian, $ G $, 
for $ \{ e_j\}_{j=1}^L $ is now
bounded by $ C^2 \max | \lambda_j | /\epsilon^2 $ so we can use $G$, 
and its inverse, to 
form a well posed Grushin problem. 
For that we change $ e_j $ to $ e_j/\sqrt{\lambda_j }$, and
denote by $ \pi $ the {\em orthogonal projection} onto the 
span of $ e_j$'s.
 We easily see the following
\begin{lem}
\label{l:3}
Condition \eqref{eq:t} implies
\[ \| ( 1 - \pi ) T \| \leq  \delta + \epsilon \| T \|  \,.\]
\end{lem}
\begin{proof}
In the notation of \eqref{eq:t} we write
\[ \tilde r = \sum_{j=M+1}^N t_j e_j + r \,, \]
and
\[ \| T u - \sum_{j=1}^M t_j e_j \| = \| \tilde r \| \leq
\delta  \| u \| + \left( \sum_{M+1}^N \lambda_j |t_j|^2 \right)^{\frac12} \,.\]
Since $ \lambda_j \leq \epsilon^2 / C^2 $, and $ \| \vec t \|_{\ell^2 } \leq    C \| T u \| $ the estimate follows.
\end{proof}
We can now  proceed as in Proposition \ref{p:3}.

\subsection{A typical estimate}
\label{sest}

Specific application of the Grushin problem scheme -- see for instance 
\S \ref{ae} -- involve estimates, often depending on a parameter. 
We would like to illustrate this
 in a situation loosely related to the
approximation scheme described in \S \ref{bfga}, and more concretely to the
example in \S \ref{aens}.

Let us assume that $ P = P ( h ) \; : \; \HH \rightarrow \HH $ is a bounded operator. 
Suppose that there exist two orthogonal projections $ \pi_\pm = \pi_\pm ( h ) 
\; : \; \HH \rightarrow \HH  $ satisfying
\begin{gather}
\label{eq:se1}
\begin{gathered}
\end{gathered}\| P^* \pi_-  \|\,,\   \| P \pi_+ \| = {\mathcal O} ( h ) \,, 
\ \ \| \pi_- P ( I - \pi_+ ) \| = o(h)\,, \ \  \| \pi_+ P^* ( I - \pi_- ) \| = o(h) \,, \\ 
\| P ( I - \pi_+ ) u \| \geq h \| ( I - \pi_+ ) u \|\,, \ \ 
\|  P^* ( I - \pi_- ) u \| \geq h \| ( I - \pi_- ) u \| \,.
\end{gather}
Here by $ a( h ) = o(h) $  we mean that $ \lim_{h\rightarrow 0 } a (h )/h  = 0 $.

We then have 
\begin{prop}
\label{p:se}
With $ P $ and $ \pi_\pm $ with the properties described 
above we define 
\begin{gather*}
  {\mathcal H}_\pm \defeq \im  \pi_\pm \,, \ \ 
R_- \; : \; {\mathcal H}_-  \hookrightarrow {\mathcal H}
\,, \ \ R_+ \; : \; \HH \stackrel{\perp}{\rightarrow} \HH_+ \,, \\
R_-^* R_- = \Id_{\HH_-} \,, \ \ R_- R_-^*  = \pi_-  \,, \ \  
R_+ R_+^* = \Id_{\HH_+} \,, \ \ R_+^* R_+  = \pi_+  \,. \end{gather*}
Then for $ h $ small enough the Grushin problem 
\begin{gather*}
\left\{ \begin{array}{ll} P u + R_- u_- & = v \\
                         R_+ u & = v_+ \end{array} \right. 
\end{gather*}
is well posed and
\begin{equation}
\label{eq:se2}
 h \| u\| + \| u_- \| \leq C ( \| v \| + h \|v _+ \| ) \,,
\end{equation}
where $ C $ is independent of $ h$.
\end{prop}
\begin{proof}
We start by rewriting our Grushin problem as 
\begin{gather}
\label{eq:se3}
\begin{gathered}
\left\{ \begin{array}{ll} P \tu + R_- u_- & = \tv \\
                         R_+ \tu & = 0 \end{array} \right. \\
\tu \defeq u - R_+^* v_+   \,,  \ \ \tv \defeq v - P R_+^* v_+ \,.
\end{gathered}
\end{gather}
We observe that $ \pi_- R_- = R_- $ and $ R_+ \pi_+ = R_+ $. 
Taking the inner product of the first equation in \eqref{eq:se3} with $ P \tu $ 
and using $ ( I - \pi_+ ) \tu = \tu $ in \eqref{eq:se1},
gives 
\begin{equation}
\label{eq:se4} \begin{split}
\frac12 \left( h^2 \| \tu \|  + \| P \tu \|^2 \right) & \leq \| P \tu \|^2   \leq  \| P \tu \| \| \tv \|  - 
\Re \langle P^* \pi_- R_- u_- , \tu \rangle  \\
& \leq \frac14 \| P \tu \|^2 +  \| \tv \|^2 + o(h) \| u_-\|_{{\mathcal H}_-}
 \| \tv \| \\ & 
  \leq 
 \frac14 \| P \tu \|^2 + 2 \| \tv \|^2 + o(h^2) \| \tu \|^2 + o(1) \| u_- \|
_{{\mathcal H}_-} ^2 \,,
\end{split} \end{equation}
Hence also
\[ \| u_- \|_{{\mathcal H}_-}  = \| R_- u_- \|
\leq  \| P \tu \| +  \| \tv \| \leq C \| \tv \| + o(1) \| u_- \|_{{\mathcal H}
_-} \,, \]
and consequently
\[ \| u _- \|_{{\mathcal H}_-} + h \| u \| \leq C \| \tv \| \,. \]
To obtain \eqref{eq:se2} we estimate  $ \| \tv \| $ using 
$ \| P R_+ ^* v_+ \| = {\mathcal O} ( h ) \| v_+ \|_{{\mathcal H}_+} $:
\[ \| \tv \| \leq \| v \| + {\mathcal O} ( h ) \| v_+ \|_{{\mathcal H}_+} \,.\]
Since by definition $ \| u \| \leq \| \tu \| + \| v_+ \|_{{\mathcal H}_+} $, the estimate
follows. 

This shows the injectivity and to see the surjectivity we apply the same 
proof to the adjoint Grushin problem, observing that the assumptions are
symmetric. 
\end{proof}

The estimate \eqref{eq:se2} is natural and appears under different 
assumptions (see for instance \cite[Lemma 5.2]{SjZw0} for another
elementary abstract estimate). In the proof we could have considered 
$ h = 1 $ since we can scale $ h $ out of the hypotheses:
\[ \left( \begin{array}{ll} h^{-1} & 0 \\
\ 0 & 1 \end{array} \right)  \left( \begin{array}{ll} P(h ) & R_-  \\
R_+  & \ 1 \end{array} \right)  \left( \begin{array}{ll} 1 &  0 \\
 0 & h \end{array} \right) =
 \left( \begin{array}{ll} h^{-1} P(h ) & R_-  \\
\ \ R_+  & \ 1 \end{array} \right) \,.\]

\subsection{Application to pseudospectral estimates}
\label{pse}

To see the estimates of \S \ref{sest} 
in use we relate them to the example in \S \ref{aens}. The general phenomenon
observed there is the growth of the resolvent of a non-normal operator 
away from the spectrum, and the consequent instability of eigenvalues -- 
see \cite{Sj},\cite{Tr},\cite{ZwPse}.

Thus consider a general $ n \times n $ matrix $ A $.
Let us then put $ P = P ( \lambda ) = A - \lambda $, and 
\[ \pi_- = \pi_- ( \lambda ) \defeq \bbbone_{ P( \lambda ) P( \lambda )^* 
\leq h^2 } \,, \ \ 
\pi_+ = \pi_+ ( \lambda ) \defeq
\bbbone_{ P( \lambda) ^*P ( \lambda )  \leq h^2 } \,, \]
where for a selfadjoint matrix $ B $, $ \bbbone_{ B \leq r } $ denotes
the orthogonal projection on the span of eigenvectors of $ B $ with eigenvalues 
less that or equal to $ r $.

We then see that the hypothesis \eqref{eq:se1} are satisfied:
\begin{gather*}
\| P \pi_+ u \|^2 = \langle P^* P \pi_+ u , \pi_+ u \rangle 
\leq h^2 \| \pi_+ u \| \,, \ \ 
\| P^* \pi_- u \|^2 = \langle P P^* \pi_- u , \pi_- u \rangle 
\leq h^2 \| \pi_- u \| \,, \\  \| P ( I - \pi_+ ) u \|^2 
= \langle P^* P ( I - \pi_+ ) u , ( I - \pi_+ )u \rangle
 \geq h^2 \| ( I - \pi_+ ) u \|\,, \\   \| P^* ( I - \pi_- ) u \|^2 
= \langle P P^* ( I - \pi_- ) u , ( I - \pi_- )u \rangle
 \geq h^2 \| ( I - \pi_- ) u \|\,, \\ 
\ \  \pi_- P ( I - \pi_+ ) = 0 \,, \ \  \pi_+ P^* ( I - \pi_- )  = 0 \,. 
\end{gather*}
To see the last identities we note that
\[ P \; : \; \ker ( P^*P - r)  \; \longrightarrow \; \ker( P P^* - r ) \,, 
\ \  P^* \; : \; \ker ( PP^* - r) \; \longrightarrow \; \ker( P^* P - r ) 
\,, \,, \]
so that 
\[ \bbbone_{ P P^* \leq h^2 } P \bbbone_{ P^*P > h^2 } = 0 \,,
\ \  \bbbone_{ P^* P \leq h^2 } P \bbbone_{ PP^* > h^2 } = 0 \,.\]
This shows that we can apply Proposition \ref{p:se} and consequently that
the Grushin problem
constructed there has the inverse:
\[  \left( \begin{array}{ll} E ( \lambda , h )  & E_+ ( \lambda ,h )   \\ 
E_-  ( \lambda ,h) 
& E_{-+} ( \lambda, h  ) 
\end{array} \right) = \left( \begin{array}{ll} {\mathcal O} ( 1/h) & {\mathcal O}(1) \\
\ \ {\mathcal O} ( 1 ) & {\mathcal O} ( h ) \end{array} \right) \; :
\CC^n \oplus \CC^{n(\lambda, h) } \; \longrightarrow \CC^n \oplus \CC^{n(\lambda, h) } \,,
\]
where $ n ( \lambda , h ) = \tr \bbbone_{ ( A - \lambda )^* ( A - \lambda ) \leq h^2 } $.
Using \eqref{eq:schur1} we see in particular that
\begin{equation}
\label{eq:pse}
 \| ( A - \lambda )^{-1} \| \simeq  \| E_{-+}(\lambda, h )^{-1} \| 
+ {\mathcal O} ( 1/h ) 
\,,\end{equation}
since $ \| E_+ v_+ \| \simeq \| v_+ \|$, and $ \| E_-^* u_- \| \simeq \| u_- \| $,
where $ a \simeq b $ means that $ b/C \leq c \leq C b $ for a constant independent of 
$ h $.

In the example presented in \S \ref{aens} where $ A $ was equal to a
Jordan block matrix, we can take any $ |\lambda|^n \ll h < |\lambda| 
 $ to obtain a Grushin problem
with $ n ( \lambda , h ) = 1 $ and  $ E_{-+} ( \lambda , h ) = \lambda^n $.

\section{Trace formul{\ae}}
\label{tf}

\subsection{Basic idea}
\label{tfbi}

Suppose that $ P = P ( z ) $. Writing $ \partial_z A ( z) = \dot A ( z ) $
we have 
\[ \dot {\mathcal E} ( z ) = - {\mathcal E} ( z) \dot {\mathcal P} ( z ) 
{\mathcal E} ( z ) \,,\]
which gives 
\begin{equation}
\label{eq:der}
E_- ( z )  \dot{P} E_+ ( z ) 
= - \dot E_{-+} ( z ) - E_- (z ) \dot R_- ( z ) E_{-+} ( z ) 
- E_{-+} ( z ) \dot R_+ ( z )  E_+ ( z ) 
\,.
\end{equation}
We recall that, formally,
\[ P(z)^{-1} = E( z ) - E_+ ( z ) E_{ -+}( z)^{-1} E_-(z )  \,.\]
Hence, assuming that we have no difficulty in taking traces, we obtain
\begin{equation}
\label{eq:der'}
\tr \dot P ( z ) P ( z )^{-1} = \tr \dot E_{-+} ( z ) E_{-+} (z)^{-1} 
+ \tr E_- ( z ) \dot R_- ( z ) + \tr \dot R_+ ( z ) E_+ ( z ) 
+ \tr \dot P E \,, 
\end{equation}
which is a special case of Lemma \ref{l:2}. 
This gives
\begin{prop}
\label{p:4}
Suppose that $ P = P ( z ) $ is a family of Fredholm operators 
depending holomorphically on $ z \in \Omega $
where $ \Omega \subset \CC $ is a connected open set. 
Suppose also that
the operators $ R_\pm = 
R_\pm ( z ) $ are of finite rank, depend holomorphically on $ z \in \Omega $,
the corresponding Grushin problem in well posed for $ z \in \Omega $,
and that  $ E_{-+} ( z_0 )^{-1} $
is invertible at some $ z_0 \in \Omega $.  Let $ g $ be  holomorphic 
 in $ \Omega$. 
Then for any curve $ \gamma $ homologous to $ 0 $ in $ \Omega $, and on 
which $ P ( z )^{-1} $ exists
$ \int_\gamma \dot P ( z ) P ( z ) ^{-1} g ( z ) dz $ is of trace class and 
we have
\begin{equation}
\label{eq:p4}
 \tr \int_\gamma  \dot P ( z ) P ( z )^{-1} g( z ) dz = 
\tr \int_\gamma \dot E_{-+} ( z ) E_{-+} (z)^{-1}  g ( z ) dz \,. 
\end{equation}
\end{prop}
\begin{proof}
Since $  E_{-+}^{-1} $ is a finite matrix for $ z \in \gamma $ we 
have that
\[ \int_\gamma \dot P ( z ) P ( z ) ^{-1} g ( z ) dz =
- \int_\gamma  \dot P ( z ) 
E_{+} ( z ) E_{-+} ( z ) ^{-1} E_{-} ( z) g ( z ) dz \,\]
is an operator of trace class, and, arguing as we did before the 
statement of the proposition we obtain \eqref{eq:p4}.
\end{proof}

The condition that $ E_{-+} ( z) $ is a finite matrix is often too 
restrictive. 
To illustrate this in a simple example we use the results of \S \ref{sebv}.
Let $ \Delta_D $ and $ \Delta_N $ be the Dirichlet and Neumann Laplacians on 
a bounded domain $ X $, with a smooth boundary $ \partial X $. 
We now put $ P_\bullet ( z ) = - \Delta_\bullet - z $, $ \bullet = D, N $. 
As described in \S \ref{sebv} we have a well posed problem for $ P_D ( z ) $
if $ P_N^{-1}( z) 
 $ exists, and in that case $ E_{-+}(z) = N(z)  $, the Neumann to 
Dirichlet operator. Similary we have 
a well posed problem for $ P_N ( z) $
if $ P_D^{-1} ( z ) $ exists, and in that case $ E_{-+}(z) = N(z)^{-1} $, 
the Dirichlet to Neumann operator. Hence if $ \gamma_D $ is a contour
homologous to $ 0 $ in the region where $ P_N ( z )^{-1}$ exists we
get $ - \tr \int_{\gamma_D} P_D ( z )^{-1} dz = \tr \int_{\gamma_D } \dot N( z)
N( z ) ^{-1} dz $. Strictly speaking we cannot apply Proposition 
\ref{p:4} directly but as $ N (z ) $ is a Fredholm operator we can locally
use an iterated problem with $ R_\pm $ of finite rank. Our contour
can be made a sum of contours involving only these local problems.

Similarly we have $ \int_{\gamma_D } P_N ( z )^{-1} dz = 0$. 
We can consider an analogous contour $ \gamma_N $ and write any contour 
$ \gamma $ as $ \gamma_D + \gamma_N $. This leads to the following 
formula:
\begin{equation}
\label{eq:dn}
 \tr \int_\gamma \left(  ( - \Delta_N - z )^{-1} - 
( - \Delta_D - z)^{-1} \right) dz  = \tr \int_\gamma N( z )^{-1} \frac{d}{dz} 
N( z ) dz \,,\end{equation}
where $ N ( z ) $ is the Neumann to Dirichlet map for $ - \Delta - z $. 
For a non-trivial application of a similar idea in the context of resonances
for the elastic Neumann problem see the work of Vodev and the first author
\cite{SjVo}.

\subsection{Classical Poisson formula}
\label{tfpf}

To present an application of Proposition \ref{p:4} we use it to  derive 
the classical Poisson
summation formula:
\begin{equation}
\label{eq:pois}
\sum_{ n \in \ZZ }  f ( n)  = \sum_{ m \in \ZZ } \hat f ( 2 \pi m ) \,, \ \ 
\ \ \hat f \in \CIc ( \RR ) \,, \ \ 
\hat f ( \xi ) \defeq \int f( x ) e^{ - i x \xi} dx \,.  \end{equation}
Our proof here might well be the most complicated derivation of \eqref{eq:pois}
but as will be indicated in \S \ref{aegt} it lends itself to far 
reaching generalizations. 

We start by rewriting \eqref{eq:pois} using the operator $ P = hD_x $ on 
$ \RR / ( 2 \pi \ZZ ) $: 
\begin{equation}
\label{eq:poiss}
 {\mathrm{tr}} \; f ( P / h ) = \frac{1 }{ 2\pi i }
\sum_{ |k|\leq N } \int_{\mathbb R} f (z/h) \left( e^{ 2 \pi i z / h }
\right)^k \frac{d}{dz} \left( e^{ 2 \pi i z / h } \right) dz \,,
\end{equation}

The left hand side there
can be written using the usual functional calculus based on Cauchy's formula:
\begin{equation}
\label{eq:2.cau}
\tr f  \left( \frac{P}{ h } \right) =  \frac{1}{ 2 \pi i}
\tr \int_\Gamma f\left ( \frac{z}{h} \right) ( P - z )^{-1} dz \,, \ 
\ \Gamma = \Gamma_+ - \Gamma_- \,, \ \ \Gamma_{\pm } = \RR \pm i R \,,
\end{equation}
where we take the positive orientation of $ \RR $ and $ R > 0 $ is an 
arbitrary constant. We make an assumption on the support of the Fourier
transform on $ f$:
\begin{equation}
\label{eq:2.f}
\supp \hat f  \subset ( - 2 \pi N , 2 \pi N ) \,.
\end{equation}

We can now use the Grushin problem \eqref{eq:2.gr} and its inverse 
given by \eqref{eq:ee}. Applying Proposition \ref{p:4} 
with $ P ( z ) = (i/h)(P - z) $ and $ g ( z ) = f ( z /h ) $ we obtain
\[ \tr f  \left( \frac{P}{ h } \right) =  - \frac{1}{ 2 \pi i}
 \int_\Gamma f\left ( \frac{z}{h} \right) \tr  
\partial_z E_{- +} ( z)  E_{ - + } ( z) ^{-1}  dz \,. \]
We now use the expression for $ E_{ - + } $ from \S \ref{fogto} to write 
\[ \begin{split} \tr f  \left( \frac{P}{ h } \right) = &  \; 
 \frac{1}{ 2 \pi i}
 \int_{\Gamma_+}  f\left ( \frac{z}{h} \right) \tr  
\partial_z M( z , h )   ( I - M ( z,h ) )  ^{-1}  dz  \\
&  + \;  
\frac{1}{ 2 \pi i } 
\int_{\Gamma_{ -} } 
  f\left ( \frac{z}{h} \right) \tr  
\partial_z M( z , h ) M ( z, h )^{-1}  ( I - M ( z,h )^{-1} )  ^{-1}  dz 
\,,   \end{split} \]
$ M ( z , h ) = \exp ( 2 \pi i z / h ) $. The assumption 
\eqref{eq:2.f} and the Paley-Wiener theorem give
\[ | \hat f ( z / h ) | \leq e^{ 2 \pi N | \Im z |/h  } \langle \Re z  / h 
\rangle^{ -\infty } \,. \]
Writing 
\[   ( I - M ( z,h ) )  ^{-1} = \sum_{ k = 0 } ^{ N -1} M ( z , h ) ^ k 
+ M ( z , h ) ^{ N} ( I - M ( z , h ))^{-1} \,, \]
for $ \Gamma_+ $, and 
\[  M ( z, h )^{-1}  ( I - M ( z,h )^{-1} )  ^{-1} =
 \sum_{ k = 1 } ^{ N} M ( z , h ) ^ {-k }
+ M ( z , h ) ^{ - N - 1} ( I - M ( z , h ))^{-1} \,, \]
for $ \Gamma_- $, we can eliminate the 
last terms by deforming the contours to imaginary infinities ($  R 
\rightarrow \infty $ in \eqref{eq:2.cau}), and this gives \eqref{eq:poiss}.

\subsection{An abstract version}
\label{tfav}

In addition to demanding finite rank of $ R_\pm $,
Proposition \ref{p:4} is restrictive in the sense that
we need to assume that the family of operators depends holomorphically on 
the parameter $ z $.
Following \cite[Appendix A]{MeSj} we present a  result without 
that assumption. 

Let $ {\mathcal H}$ be a complex Hilbert space and let us denote by
$ {\mathcal L} ( {\mathcal H }, {\mathcal H} ) $
bounded operators on $ {\mathcal H} $. 
We consider  $ \SP^1 \ni t \mapsto A ( t ) \in {\mathcal L}
( {\mathcal H}, {\mathcal H} ) $,  a $ {\mathcal C}^1 $ 
closed curve of operators, in the sense of $ A ( t) $ is strongly 
differentiable with a continuous derivative $ \SP^1 \ni t \mapsto 
\dot A ( t ) \in {\mathcal L} ( {\mathcal H }, {\mathcal H} )$.
We write $ dA = \dot A dt $, and for another such $ t \mapsto B ( t ) $,
\[ \int_{\SP^1}  B ( t) \dot A ( t ) dt = \int B dA \in {\mathcal L}
(\HH, \HH) \,. \]

If the values of $ A ( t ) $ are
taken in an open subset $ V $ of $ {\mathcal L} ({\mathcal H}, 
{\mathcal H} )  $, we will 
will say that $ A ( t ) $ is {\em contractible} in $ V $, if 
$ \SP^1 \ni t \mapsto A ( t ) \in V $ has a $ {\mathcal C}^1 $ 
extension $ \DP \ni z \mapsto A ( z ) \in V $, $ \DP = \{ |z | < 1 \} $,
$ \partial \DP = \SP^1 $. 

With this terminology we have
\begin{prop}
\label{p:5}
Suppose that 
\begin{equation}
\label{eq:p.5'}
 \SP^1 \ni t \mapsto {\mathcal P}( t ) =
\left( \begin{array}{ll} P ( t ) & R_- ( t ) \\
R_+ ( t ) & \ \ 0 \end{array} \right) \; : \; 
{\mathcal H} \oplus
 \CC^N \; \longrightarrow \;  
{\mathcal H} \oplus \CC^N\,,\end{equation}
is contractible in the set of invertible operators on $ 
{\mathcal H} \oplus \CC^N $, with  
$ \DP \ni z \mapsto \dot {\mathcal P} ( z )  $  continuous
with values in operators of trace class.

If $ P(t)  ^{-1} $ exists for all $   t \in \SP^1 $ then 
\begin{equation}
\label{eq:p.5}
\tr \int P^{-1} dP  = \tr \int E_{-+}^{-1} dE_{-+} \,,
\end{equation}
where we use the standard Grushin problem notation for the inverse of 
\eqref{eq:p.5'}.
\end{prop}
\begin{proof}
Schur's formula \eqref{eq:schur1} and the fact that $ dP $ is of trace class 
give
\[ \tr P^{-1} dP = \tr E dP - \tr E_+ E_{-+}^{-1} E_- dP \,.\]
Using \eqref{eq:der} and the cyclicity of the trace we see that
\[ \tr P^{-1} dP = \tr E_{-+}^{-1} dE_{-+} + \omega \,, \ \ 
\omega \defeq \tr ( dR_+ E_+ ) + \tr ( E_- dR_- ) + \tr ( E dP ) \,.\]
We will obtain \eqref{eq:p.5} when we show that $ \omega $, which is 
defined on the circle, extends to a closed form in the unit disc. 
Since $ t \mapsto {\mathcal P }(t)  $ is assumed to be contractible
we can use the same notation for its extension to the unit disc:
\[  {\mathcal P}( z ) =
\left( \begin{array}{ll} P ( z ) & R_- ( z ) \\
R_- ( z ) & R_{+-} ( z )  \end{array} \right) \,.\]
Lemma \ref{l:2} then shows that $ \omega $ is a restriction to 
the unit circle of a one form defined in the unit disc:
\[ \omega = \tr ( dR_+ E_+ ) + \tr ( E_- dR_- ) + \tr ( E dP ) 
+ \tr ( dR_{+-} E_{-+} ) \,.\]
To compute $ d \omega $ we note that  $ d {\mathcal E} = - {\mathcal E}
d {\mathcal P} {\mathcal E } $ and consequently
\[ \begin{split}
& - d E = E dP E + E_+ dR_+ E + E dR_- E_- + E_+ dR_{+-} E_- \,, \\
& - d E+ = E dP E_+ + E_+ dR_+ E_+ + E dR_- E_{-+} + E_+ dR_{+-} E_{-+} \,,\\
& - d E_- = E_- dP E + E_{-+ } dR_+ E + E_- dR_- E_- + E_{-+} dR_{+-} E_- \,, \\
& - d E_{-+} = E_{-} dP E_+ + E_{-+} dR_+ E_+ + E_- dR_- E_{-+} + E_{-+ }dR_{+-} E_{-+} 
\,.
\end{split} \]
Hence, using the natural notation for operator valued differential forms, we 
obtain
\[ \begin{split} 
d \omega = & \ \tr dR_+  \wedge  E dP E_+ + \tr d R_+ \wedge E_+ dR_+ E_+ \\ 
& + \; \tr  dR_+ \wedge E dR_- E_{-+} + \tr  dR_+ \wedge E_+ dR_{+-} E_{-+}  \\ 
& - \; \tr   E_- dP E \wedge  dR_-  - \tr E_{-+ } dR_+ E \wedge  dR_- \\ 
& - \; \tr E_- dR_- E_- \wedge  dR_- -  \tr E_{-+} dR_{+-} E_- \wedge  dR_- \\ 
& -   \; \tr  E dP E \wedge  dP   - \tr E_+ dR_+ E \wedge  dP  \\ 
& - \; \tr E dR_- E_- \wedge  dP  - \tr E_+ dR_{+-} E_- \wedge  dP \\ 
& + \; \tr dR_{+-} \wedge E_{-} dP E_+ + \tr dR_{+-} \wedge E_{-+} dR_- E_+ \\ 
& + \; \tr dR_{+-} \wedge  E_- dR_- E_+ + \tr dR_{+-} \wedge E_{-+} dR_{+-} E_{-+}  \, .  
\end{split} \]
For a differential 1-form, $ \mu $,  with trace class we clearly have 
$ \tr \mu \wedge \mu =0 $. That shows that in the expression for $ d\omega $
the 2nd, 7th, 9th, and 16th terms vanish. Cyclicity of the trace also shows that
the  terms in pairs: (1st,10th), (3rd,6th), (4th,14th), (12th,13th), (8th,15th)
cancel each other.
Finally,  we also have $ \tr \mu_1 \wedge 
\mu_2 = - \mu_2 \wedge \mu_1 $, and it follows that the
5th and 11th  terms  cancel each other. Thus $ d \omega = 0 $ completing the proof.
\end{proof}

It is quite possible that Proposition \ref{p:5} follows from some general
topogical facts. It is not clear what are the weakest assumptions on $ P $ and
$ d P $ to guarantee that $ \int P dP $ is of trace class. For a discussion of
one case of a weaker assumption see \cite[Appendix A]{MeSj}.

\section{Advanced examples}
\label{ae}

\subsection{Around Lidskii's perturbation theory for matrices}
\label{lr} 

In \S \ref{aens} the equation for the eigenvalues of the perturbation 
of one Jordan block is easily derived from \eqref{eq:empe}:
\[ \lambda^n - \epsilon Q_{n1} + \epsilon {\mathcal O}(\lambda ) 
+ {\mathcal O} ( \epsilon^2 ) = 0 \,, \ \ 
|\lambda | < 1 \,.\]
and the solutions are 
\[ \lambda_\ell = \epsilon^{1/n} |Q_{n1}|^{1/n}
 e^{ (2 \pi i \ell + \arg{Q_{n1}} )/n  } + 
o( \epsilon^{1/n} ) \,, \ \ 1 \leq \ell \leq n \,.\]
Here we consider $ n $ fixed and are interested in the $ \epsilon 
\rightarrow 0 $ asymptotics.

In this section we will show how the Grushin problem approach applies to 
the study of perturbation of matrices with arbitrary Jordan 
structure. We restrict ourselves to an example suggested 
by Michael Overton which according to him
contains the essential elements of the general problem studied
in \cite{Lid} and \cite{MBO}. 

Let $ J_\ell $ be the $ \ell \times \ell$ upper triangular Jordan bloc
matrix. We then consider
\begin{equation}
\label{eq:A}  
A = J_n \oplus J_n \oplus J_k \; : \; \CC^n \oplus \CC^n 
\oplus \CC^k  \; 
\longrightarrow \CC^n \oplus \CC^n 
\oplus \CC^k  \,, \ \  k < n \,, \end{equation}
that is 
\[ A = \left(\begin{array}{lll}  J_n & 0_{nn} & 0_{nk} \\
0_{nn} & \ J_n & 0_{nk} \\
0_{kn} & 0_{kn} & J_k \end{array} \right)\,,\]
where $ 0_{\ell p } $ denotes the $ \ell \times p $ zero matrix.

The Grushin problem for $ A $ is a straightforward modification of the
one for $ J_n $ in \S \ref{aens}:
\begin{gather*} R_- \; : \; \CC^3 \; \longrightarrow \CC^n \oplus \CC^n 
\oplus \CC^k  \,, \ \ 
R_+ \; : \; \CC^n \oplus \CC^n 
\oplus \CC^k  \; \longrightarrow \CC^3 \,. \end{gather*} 
We then obtain the effective Hamiltonian, $ E_{-+} (\lambda ) $ for 
$ A - \lambda $:
\[ E_{-+} ( \lambda ) = \left( \begin{array}{lll} 
\lambda^n & \ 0 & \ 0 \\
\ 0 & \lambda^n & \ 0 \\
\ 0 & \ 0 & \lambda^{k} \end{array} \right)
\,,\]
and $ E_\pm ( \lambda ) $ are similarly constructed from the three 
$ e_\pm( \lambda ) $ vectors. 

Suppose we now consider a perturbation of $ A $: 
\begin{gather} 
\label{eq:Ae}
\begin{gathered}
A_\epsilon = A + \epsilon Q \,, \ \
Q = \left( \begin{array}{lll} Q^{11} & Q^{12} & Q^{13} \\
Q^{21} & Q^{22} & Q^{23} \\
Q^{31} & Q^{32} & Q^{33} \end{array} \right) \,, \\
Q^{ij} \; : \; \CC^{n_i} \; \longrightarrow \; \CC^{n_j} \,, \ \
n_i = n\,, \ i=1,2\,, \ \ n_3 = k \,.
\end{gathered}
\end{gather}

As in \S \ref{aens} we see that the effective Hamiltonian for the perturbation
is 
\[   E_{-+}^\epsilon ( \lambda ) = 
E_{-+} ( \lambda ) - 
 \epsilon E_- (\lambda ) Q  E_+ ( \lambda )  + {\mathcal O} ( \epsilon^2) 
  \,.
\]
The effective first order perturbation is easily checked to be
\[  E_- (\lambda ) Q  E_+ ( \lambda ) 
=
 \left( \begin{array}{lll} Q_{n1}^{11} & Q_{n1}^{12} & Q_{n1}^{13} \\
\ & \ & \ \\
Q_{n1}^{21} & Q_{n1}^{22} & Q_{n1}^{23} \\
\ & \ & \ \\
Q_{k1}^{31} & Q_{k1}^{32} & Q_{k1}^{33}  
\end{array} \right)  + {\mathcal O} ( \lambda ) 
\,,\] 
where $ Q_{ij}^{pq} $ denotes the $ ij$'th entry of the matrix 
$ Q^{pq} $. 

Suppose that the matrix $ (Q^{ij}_{n1} )_{1\leq i,j, \leq 2 } $ is 
diagonalizable with eigenvalues $ q_1 $ and $ q_2 $. Then the
eigenvalues of $ A_\epsilon $ are given by the values of $ \lambda $
for which the following matrix is not invertible:
\[  \left( \begin{array}{lll} \lambda^n - \epsilon q_1 & 
\ \ \ 0  &  \ \ \epsilon \widetilde Q_{13}  \\
\ \ \ 0  & \lambda^n - \epsilon q_2 & 
\ \  \epsilon \widetilde Q_{23}  \\
\ \ \epsilon \widetilde Q_{31} &  \ \ \epsilon \widetilde Q_{32} &  
\lambda^k - \epsilon \widetilde Q_{33} \end{array}\right) + 
\epsilon {\mathcal O} ( \lambda ) + {\mathcal O} ( \epsilon^2 ) \,.
\]
Since $ k < n $, and both $k $ and $ n $  are fixed, perturbation theory gives 
\begin{prop}
\label{p:l}
The largest modulus eigenvalues of $ A_\epsilon $ for $ \epsilon $ small 
are given by 
\[  \lambda_\ell^j = \epsilon^{1/n} |q_j|^{1/n}
 e^{ (2 \pi i \ell + \arg{q_j} )/n  } + 
o( \epsilon^{1/n} ) \,, \ \ 1 \leq \ell  \leq n \,, \ \ j = 1, 2 \,, \]
where $ q_{j} $ are the eigenvalues (assumed to be distinct) of the
$ ( Q^{ij}_{n1} )_{1 \leq i,j \leq 2} $ part of the perturbation matrix
in \eqref{eq:Ae}.
\end{prop}
A finer perturbation theory for matrices of size given by the
number of distinct Jordan blocks will (most likely) give the
general results of \cite{Lid} and \cite{MBO}.

\subsection{Generalized Gutzwiller trace formula}
\label{aegt}

Trace formul{\ae} provide one of the most elegant descriptions of the
classical-quantum correspondence. One side of a formula is given 
by a trace of a quantum object, typically derived from a quantum
Hamiltonian, and the other side is described in terms of closed 
orbits of the corresponding classical Hamiltonian. 

Here we follow \cite{SjZw} 
and outline the structure of 
a formula which is derived using a {\em formal} Grushin problem. 
It is an intermediate trace formula in which
the original trace is expressed in terms of traces of quantum monodromy
operators directly related to the classical dynamics. The usual trace
formul{\ae} follow and in addition this approach allows handling effective
Hamiltonians, such as the one described in \S \ref{aeps} below. 

Let $ P $ be a semi-classical, self-andjoint, principal type operator,
elliptic in the classical sense, with symbol $ p$, and a compact 
characteristic variety, $ p^{-1} ( 0 ) $.
Let $ \gamma \subset p^{-1} (0 ) $ 
be a closed primitive orbit of the Hamilton flow of $ p$. The simplest example,
and one discussed in \S \ref{tfpf}, $ P = hD_x$, on the circle,
$ p^{-1}(0) = \{ ( x , 0 ) \; x \in \SP^1 \} \subset T^* \SP^1 $, and the 
Hamilton vector field is $ \partial_x $. More interesting examples are
$ P = - h^2 \Delta_g - 1 $ on a compact Riemannian manifold, or 
$ P = -h^2 \Delta + V ( x ) $ with a suitable $ V $ on $ \RR^n $.

We can define the {\em monodromy operator}, $ M ( z, h ) $ for $ P - z $ 
along $ \gamma$, acting on functions in one dimension lower, that is, 
on functions on the transversal 
to $ \gamma $ in the base. We then have 
\begin{thm}
\label{t:0}
Suppose that there exists a neighbourhood of $ \gamma $, 
$ \Omega $, satisfying the condition
\begin{equation}
\label{eq:cond}
 m \in \Omega  \ \text{and} \ \exp t H_p ( m ) = m \,, \ 
p ( m ) = 0 \,, \ \ 0 <  |t| \leq T N 
\ \Longrightarrow \ m \in \gamma \,, \end{equation}
where $ T $ is the primitive period of $ \gamma$. 
If $ \hat f \in 
{\mathcal C}^\infty_{\rm{c}} ( {\mathbb R} )$, 
$ {\rm supp}\; \hat f \subset ( - NT  , NT  ) \setminus \{ 0\}$, 
$ \chi \in \CIc (\RR ) $, and $ A \in \Psi_h^{0,0} ( X ) $ is a microlocal 
cut-off to a sufficiently small neighbourhood of $ \gamma $, then
\begin{equation}
\label{eq:main}
 {\mathrm{tr}} \; f ( P / h ) 
\chi ( P ) 
A = \frac{1 }{ 2\pi i } 
\sum_{  -N - 1}^{N - 1 } {\mathrm {tr}} \; 
\int_{\mathbb R} f (z/h) M(z, h)^k \frac{d}{dz} M ( z , h )
 \chi ( z) 
dz + {\mathcal O} ( h^\infty ) \,, \end{equation}
where $ M ( z, h )$ is the semi-classical monodromy operator associated to 
$ \gamma $. 
\end{thm}
The dynamical assumption on the operator means
that in a neighbourhood of $ \gamma $ there are no other closed orbits of 
period less than $ T N$, on the energy surface $ p  = 0 $. We avoid 
a neighbourhood of $ 0 $ in the support of $ \hat f $ to avoid the
dependence on the microlocal cut-off $ A $.

The monodromy operator quantizes the Poincar{\'e} map for $ \gamma $ and 
its geometric analysis gives the now standard trace formul{\ae} of 
Selberg, Gutzwiller and Duistermaat-Guillemin.
The term $ k = -1 $ corresponds
to the contributions from ``not moving at all'' and the other terms to 
contributions from going $ |k+1| $ times around $ \gamma $, in the positive
direction when $ k \geq 0$, and in the negative direction, when $ k < -1 $.
For non-degenerate orbits the analysis of  the traces on monodromy operators 
recovers the usual semi-classical trace formul{\ae} in 
our general setting -- see \cite[Theorem 3]{SjZw}.

The proof of the formula follows the lines of the proof of classical 
Poisson formula presented in \S \ref{tfpf}.
In the general situation where the circle is replaced by a closed 
trajectory of a real principal type operator 
we can proceed similarly but now {\em microlocally}
in a neighbourhood of that closed orbit on an energy surfarce.
The contour integral formula \eqref{eq:2.cau} is replaced by the 
Dynkin-Droste-Helffer-Sj\"ostrand formula (see \cite[Chapter 8]{DiSj})
\begin{equation}
\label{eq:2.cau'}
\tr f  \left( \frac{P}{ h } \right) \chi ( P ) A = -  \frac{1}{ \pi }
 \int_{\CC}  f\left ( \frac{z}{h} \right) \bar \partial_z
\tilde \chi ( z) ( P - z )^{-1} A\;  {\mathcal L}( dz)  \,, 
\end{equation}
where $ \tilde \chi $ is an {\em almost analytic extension} of
$ \chi $, that is an extension satisfying $ \bar \partial_z \chi ( z) 
= {\mathcal O} ( | \Im z |^{\infty } ) $ -- see \cite[Sect.6]{SjZw}
and we want to proceed with a similar reduction to the effective 
Hamiltonian given in terms of a suitably defined {\em monodromy operator}.

To construct the monodromy operator we fix two different points on $ \gamma$,
$ m _0 $, $ m_1$ (corresponding to $ 0 $ and $ \pi $ in
\eqref{eq:2.gr}-\eqref{eq:2.mon}),
and their disjoint neighbourhoods, $ W_+ $ and $ W_- $ respectively.
We then consider local kernels of $ P - z$ near $ m_0 $ and $ m_1 $ (that
is, sets of distributions satisfying $ ( P -z ) u = 0 $ near $ m_i $'s),
$ \ker_{m_j} ( P - z ) $, $ j = 0, 1$,  with elements microlocally 
defined in $ W_\pm $.
and the forward and backward solutions:
\[ I_\pm ( z) \; : \; \ker_{ m_0 } ( P - z) \ \longrightarrow \
\ker_{ m_1 } ( P - z ) \,.  \]
We then define the {\em quantum monodromy operator}, $ {\mathcal M} (z)$
by
\renewcommand\thefootnote{\ddag}%
\begin{gather*}
 I _{-}  ( z ) {\mathcal M} ( z ) = I_+ ( z ) \,, \ \ \ \
{\mathcal M} ( z) \; : \;   \ker_{ m_0 } ( P - z)
\ \longrightarrow \  \ker_{ m_0 } ( P - z) \,. \end{gather*}
The operator $ P $ is assumed to be self-adjoint with respect to some
inner product $ \langle \bullet, \bullet \rangle$, and we define
the {\em quantum flux} norm on $ \ker_{m_0} ( P - z) $ as 
follows\footnote{See \cite{HeSj2} for an earlier mathematical 
development of this basic quantum mechanical idea.}: let $ \chi $ be 
a microlocal cut-off function, with basic properties of the function
$ \chi $ in the example. Roughly speaking $ \chi $ should supported near
$ \gamma $ and be equal to one near the part of $\gamma $ between $ W_+$ and 
$ W_- $. We denote by $ [ P , \chi ]_{W_+} $ the part of the
commutator supported in $ W_+$, and  put
\[  \langle u , v \rangle_{\rm{QF}}  \stackrel{\rm{def}}{=} 
\langle [ ( h / i ) P , \chi ]_{W_+} u , v \rangle \,, \ \ u ,v 
\in \ker_{m_0} ( P - z ) \,. \]
As can be easily seen this norm is independent of the choice of 
$ \chi $.  This independence leads
to the unitarity of $ {\mathcal M} ( z ) $:
\[  \langle {\mathcal M} ( z)  u ,  {\mathcal M} ( z )  u \rangle_{\rm{QF}}
=  \langle u , u  \rangle_{\rm{QF}} \,, \ \ u \in \ker_{m_0} ( P - z )  
\,. \]
For practical reasons we identify $ \ker_{ m_0 } ( P - z) $ with 
$ {\mathcal D}' ( \RR^{n-1} ) $, microlocally near $ ( 0 , 0 )$, and
choose the idenfification so that the corresponding monodromy map is
unitary (microlocally near $ ( 0 , 0 ) $ where $ ( 0, 0 )$ corresponds
to the closed orbit intersecting a transversal identified with 
$ T^* \RR^{n-1} $). This gives
\[ M ( z , h ) \; : \; {\mathcal D} ' ( \RR^{n-1} ) \ \longrightarrow
\ {\mathcal D}' ( \RR^{n-1} ) \,, \]
microlocally defined near $ ( 0, 0 ) $ and unitary there.
This is the operator appearing in Theorem \ref{t:0} and it shares
many properties with its simple version $ \exp(2 \pi i z / h ) $ appearing in 
\eqref{eq:2.mon}
for $ {\mathbb S}^1$. 

\renewcommand\thefootnote{*}%

As in \S \ref{fogto} we can construct a Grushin problem with the 
effective Hamiltonian given by $ E_{-+} ( z , h) = I - M ( z, h ) $.
However, now the problem is {\em formal}, that is all the inversion formul{\ae}
are only valid microlocally\footnote{For a review of this
important notion see \cite[Section 3]{SjZw}. Roughly speaking it corresponds
to a localization of 
the behaviour of quantum states to relevant subsets of classical 
phase space. It does not guarantee global well-posedness in an {\em honest}
Hilbert space sense.} near $ \gamma $. Since in Theorem \ref{t:0} we are
interested in taking traces, and not, for instance, locating eigenvalues
or resonances, that is sufficient. 

Nevertheless, 
as one striking application of this point of view we can explain the way in 
which complex quasi-modes manifest themselves on compact manifolds \cite{ISZ},
a phenomenon which was already explicitely or implicitely noted in the works of
Paul-Uribe, Guillemin, and Zelditch -- see \cite{Zel} and references given there.

To explain it, let us recall the now classical fact (Lazutkin, Ralston, 
Colin de Verdi\`ere, Popov)
that for an elliptic closed geodesic on a compact
manifold $ M $ one can construct approximate eigenfunctions concentrating
on that trajectory, and that the corresponding 
approximate eigenvalues are close to actual
eigenvalues with arbitrary polynomial accuracy as energy increases. 
When the trajectory is hyperbolic that procedure no longer makes sense
as the formal construction of quasi-modes gives complex numbers. That 
can lead to the construction of resonances in scattering situations
(Ikawa, G\'erard, Sj\"ostrand-G\'erard) but cannot have a direct 
spectral interpretation when the manifold is compact. Despite that they 
make a direct appearance when traces are considered and we have 
the following consequence of recent work on inverse spectral problems 
(see \cite{Zel} and \cite{ISZ}):
\begin{thm}
\label{t:q}
Let $ M $ be a compact Riemannian Riemannian manifold and $ \gamma $ a closed
hyperbolic trajectory of primitive length $ L_\gamma $. 
Let $ \lambda_j^2 $ denote the sequence of eigenvalues
of the Riemann-Beltrami operator, 
$ \mu_{k} $ the sequence of complex 
quasi-modes associated to the trajectory $ \gamma $, $ 0 < \Im \mu_{k} $
(well defined modulo $ {\mathcal O} ( |\Re \mu_k|^{-\infty } ) $).
Suppose that for any $ m \in \ZZ \setminus \{ 0 \}$, 
$ m L_\gamma $ is different from the length of any  closed geodesic on  $ M $
which is not an iterate of $ \gamma$.
Then, for any $ m \in \ZZ \setminus \{ 0 \} $ there exists a neighbourhood $ U_m $
of $ m L_\gamma $ such that 
\[  \sum_j e^{i \lambda_j t } - \sum_{k } e^{ i \mu_k t } \in \CI( U_m ) \,,\]
where both sums are meant in the sense of distributions on $ \RR $, and 
$  \sum_{k } e^{ i \mu_k t } $ is defined only modulo 
$ \CI ( \RR \setminus 0 ) $.
\end{thm}

In our approach, epecially in view of Grushin reductions to the
effective Hamiltonians, it is important that we can consider operators
with non-linear dependence on the spectral parameter. In that case, motivated
by Proposition \ref{p:4}, 
the left hand side of \eqref{eq:main} is replaced by 
\[ \frac1\pi  {\mathrm{tr}} \;  
\int f ( z / h )
 \bar \partial_z \left[ \tilde \chi ( z )\; 
\partial_z P (z) \;  P ( z) ^{-1} \right] A
{\mathcal L} ( d z ) \,,\]
which for $ P ( z) = P - z$ reduces to \eqref{eq:2.cau'}. For a generalized
version we refer to \cite[Theorem 2]{SjZw}.

Finally we point out that 
the semi-classical Grushin problem 
point of view 
taken here, when translated to 
the special case of $\CI$-singularities/high energy regime,
is close to that of Marvizi-Melrose and Popov (see references in \cite{SjZw})
In those works the trace of the wave group was reduced to the study
of a trace of an operator quantizing the Poincar\'e map.

\subsection{Peierls substitution}
\label{aeps}
In this section we will follow \cite{HeSj1} 
to show how the Grushin reduction leads to a
natural mathematical explanation of the celebrated Peierls substitution from 
solid state physics. It gives an effective Hamiltonian for a crystal in 
a magnetic field. For simplicity of the presentation 
we will consider the case of dimension 
two only, and of the first spectral band -- we refer to \cite{HeSj1} and
\cite{HeSj2} for
the general case and for references to the vast literature on the subject.

First we need to consider the case of no magnetic field. Mathematically 
this corresponds to considering a Schr\"odinger operator with a 
periodic potential:
\[ P_0 = -\Delta + V \,,  \ \ V \in \CI ( \RR^2 ) \,, \ \
V ( x + \gamma ) = V ( x ) \,, \ \ \gamma \in \Gamma \,, \]
where $ \Gamma $ is a lattice in $ \RR^2 $. In other words,
\begin{equation}
\label{eq:tr}
T_\alpha P_0 = P_0 T_\alpha \,, \ \  T_\alpha u ( x ) \defeq u ( x - \alpha ) \,.
\end{equation}
The operator $ P $ is unitarily equivalent to a direct integral of {\em Floquet
operators}, $ P_\theta $, acting as $ P $ on $ {\mathcal H}_\theta $:
\[ {\mathcal H}_{\theta} \defeq \{ u \in L^2_{\rm{loc} } ( \RR^2 ) \; : \; 
\forall \gamma \in \Gamma \ u ( x - \gamma ) = e^{ i \langle \theta, \gamma \rangle }
u ( x ) \} \,, \ \ \theta \in \RR / \Gamma^* \,, \]
where $ \Gamma^* $ is the dual lattice of $ \Gamma $: $ \gamma^* \in \Gamma^* 
\Longleftrightarrow \langle \gamma^*, \alpha \rangle \in 2 \pi \ZZ $ for all 
$ \alpha \in \Gamma $. We denote by $ E $ and $ E^*$ the fundamental domains of
$ \Gamma $ and $ \Gamma^* $ respectively. Explicitely, 
\begin{gather*}
{\mathcal B} P_0 {\mathcal C} = \int^{\oplus} P_\theta d \theta \,, \\
({\mathcal B } f )( x, \theta ) = \sum_{\gamma \in \Gamma } 
e^{ - i \langle \theta, \gamma \rangle } f ( x - \gamma ) \,, \ \ 
( {\mathcal C} g ) ( x ) = \frac{1}{\vol( E^* ) } \int_{E^* } g ( x , \theta ) 
d \theta \,, \\
{\mathcal B } \; : \; L^2 ( \RR^2 ) \; \longrightarrow L^2 ( \RR^2/\Gamma^* , 
{\mathcal H}_\theta ) \,, \ \ {\mathcal C} = {\mathcal B}^* = {\mathcal B}^{-1}\,.
\end{gather*}
The spectrum of $ P_0 $ is absolutely continuous and equal to $ \bigcup_{k \in 
\NN } \{ \lambda_k ( \theta ) \; : \; \theta \in \RR^2 / \Gamma^* \} $, where
$ \{ \lambda_k ( \theta )\}_{k=1}^\infty $ is the sequence of eigenvalues of 
$P_\theta $. Each interval in the union is referred to as a band and we assume
that the first band is disjoint from all the other bands.

We now want to find a Grushin problem for $ P - z $ which will be well posed
near the first band. It turns out (see \cite[Lemma 1.1]{HeSj})
that one can choose 
$ \phi ( x , \theta ) $, $ P_\theta \phi ( x , \theta ) = \lambda_1 ( \theta ) 
\phi ( x , \theta ) $, to be holomorphic, as a function of  $ \theta $, 
in a complex neighbourhood of $ \RR^n /\Gamma^* $. That implies that 
\[ \phi_0 ( x ) \defeq ( {\mathcal C} \phi ) ( x ) \,, \]
has very nice properties: $ | \partial_x^\alpha \phi_0 ( x ) | 
\leq C_\alpha e^{ - |x|/ C} $.
We then define
the following Grushin problem:
\begin{gather}
\label{eq:pgr1}
\begin{gathered}
{\mathcal P}_0 ( z ) = 
\left( \begin{array}{ll} P_0 - z & R_-^0 \\
\ \ R_0^+ & \ \ 0 \end{array} \right) \; : \; H^2 ( \RR^2 ) \oplus
\ell^2 ( \Gamma ) \; \longrightarrow L^2 ( \RR^2 ) \oplus \ell^2 ( \Gamma ) \,,\\
R_+^0 \defeq (R_-^0)^* \,, \ \  R_-^0 u_- ( x ) \defeq \sum_{\gamma \in \Gamma }
u_- ( \gamma ) \phi_0 ( x  - \gamma ) \,. 
\end{gathered}
\end{gather}
It is not hard to see that this problem is well posed for $ z $ close to the
first band and away from all the other bands. The effective Hamiltonian is
given by 
\[ ( E_{-+}^0( z) v_+ ) ( \alpha ) = \sum_{\beta \in \Gamma}
( z \delta_{\alpha, \beta} - \widehat E ( \alpha - \beta ) ) v_+ ( \beta) \,, 
\ \ \widehat E ( \gamma ) = \frac1{\vol(E^*)} \int_{E^* } \lambda_1 ( \theta ) e^{ i 
\langle \theta, \gamma \rangle } d \theta \,,\]
which is unitarily equivalent to the multiplication by $ z - \lambda_1 ( \theta ) $,
the obvious effective Hamiltonian near the first band.

The Grushin problem \eqref{eq:pgr1} does have the advantage of being
stable under small perturbations and we will see it when the magnetic field
is turned on. That correponds to adding a magnetic potential to our operator.
Here we consider only a constant weak magnetic field $ B = h dx_1 \wedge dx_2 $:
\begin{equation}
\label{eq:cmf}
 P_B = ( D_{x_1} - h x_2 )^2 + D_{x_2}^2  + V ( x) \,, \ \
D_{x_j} = \frac{1}{i} \partial_{x_j } \,, \ \  B = h dx_1 \wedge dx_2 \,. 
\end{equation}
Although the operator $ P_B $ is no longer periodic in the sense of \eqref{eq:tr}
it commutes with magnetic translations:
\begin{equation}
\label{eq:trB}
T_\alpha^B P_B = P_B T_\alpha^B \,, \ \  T^B_\alpha u ( x ) \defeq 
e^{ \frac{i}2 \langle B, x \wedge \alpha \rangle } u ( x - \alpha ) \,, \ \ 
T_\alpha ^B T_\beta^B = e^{ - i \langle B , \alpha \wedge \beta \rangle } 
T_\beta^B T_\alpha^B \,. 
\end{equation}
We now use the magnetic translations to modify 
the Grushin problem \eqref{eq:pgr1}:
\begin{gather}
\label{eq:pgr2}
\begin{gathered}
{\mathcal P}_B ( z ) = 
\left( \begin{array}{ll} P_B - z & R_-^B \\
\ \ R_B^+ & \ \ 0 \end{array} \right) \; : \; H^2_B ( \RR^2 ) \oplus
\ell^2 ( \Gamma ) \; \longrightarrow \; L^2 ( \RR^2 ) \oplus \ell^2 ( \Gamma ) \,,
\\  (R_-^B u_-)(x)   \defeq \sum_{\gamma \in \Gamma }
u_- ( \gamma ) T_\gamma^B \phi_0 (x)  \,, \ \ 
(R_+^B u) (\gamma ) 
\defeq \langle u , T_\gamma^B \phi_0 \rangle_{L^2 ( \RR^2 ) } \,, \\
H^2_B \defeq \{ u \in L^2 ( \RR^2 ) \; : \;  P_B u \in L^2 ( \RR^2 ) \} \,. 
\end{gathered}
\end{gather}
The operator $ {\mathcal P}_B ( z ) $ commutes with 
\[ \left( \begin{array}{ll} T_\gamma^B & 0 \\
\ \ 0 & \tau_\gamma^B \end{array} \right) \,, \ \  \tau_\gamma^B u (\alpha ) 
\defeq
e^{ \frac{i}{2} \langle B, \alpha \wedge \gamma \rangle } u ( \alpha - \gamma ) \,.
\]
It is shown in \cite[Proposition 3.1]{HeSj1}
that when $ h $ is small ($ B = h dx_1  \wedge dx_2 $) then $ {\mathcal P}_B ( z ) $
is invertible for $ z $ near the first band for $ P_0 $. Although it requires
some technical work, roughly speaking it follows from the invertibility of
$ {\mathcal P}_0 $ and the smallness of the magnetic field.

The inverse has the same symmetries as $ {\mathcal P}_B ( z )$ and in 
particular $ \tau_\alpha^B E_{-+} ( B, z ) = E_{-+} ( B, z ) \tau_\alpha^B $
for all $ \alpha \in \Gamma $. That implies that $ E_{-+} ( z , B ) $ 
is given by a ``twisted convolution'':
\begin{equation}
\label{eq:twist}
\left( E_{-+} ( z , B ) v_+ \right)( \alpha ) = \sum_{\beta \in \Gamma } 
e^{\frac{i}{2} \langle B , \alpha \wedge \beta \rangle } f_{B,z} ( \alpha -
\beta ) v_+ ( \beta ) \,, \ \ 
|f_{ B, z } ( \gamma ) | \leq C e^{ | \gamma | / C } \,. \end{equation}
Operators with kernels satisfying these properties form an algebra sometimes
called the algebra of {\em magnetic matrices}. In \cite[Proposition 5.1]{HeSj1}
it is shown that the inveribility of a magnetic matrix as an operator
on $ \ell^2 ( \Gamma ) $ is equivalent to its invertibility in the algebra
of magnetic matrices. Let $ {\mathcal M}_B ( f ) $ denote the magnetic 
matrix associated to an exponentially decaying function on $ \Gamma $, $ f $:
\[   {\mathcal M}_B ( f ) ( \alpha, \beta ) = 
e^{\frac{i}{2} \langle B , \alpha \wedge \beta \rangle } f ( \alpha -
\beta )  = 
e^{\frac{i}{2} h \sigma ( \alpha, \beta ) } f ( \alpha -
\beta )
\,, \]
where $ \sigma $ is the standard symplectic form on $ \RR^2 $. It is easy to check 
that 
\begin{equation}
\label{eq:mag}
   {\mathcal M}_B ( f ) \circ   {\mathcal M}_B ( g )  =
  {\mathcal M}_B ( f \; \#_B \;  g  ) \,, \ \
f\; \#_B \;  g (\gamma)  = \sum_{\alpha +\beta = \gamma } e^{ \frac{i}{2}
\langle B , \alpha \wedge \beta \rangle } f( \alpha) g ( \beta ) \,. \end{equation}

We are now getting close to the Peierls substitution which provides an 
elegant microlocal description of $ E_{-+}(z, B ) $. We can take
the Fourier transform of an exponentially decaying function on $ \Gamma $, $ f$,
\[ \widehat f ( \theta ) \defeq \sum_{\gamma\in \Gamma } e^{i \langle \theta, \gamma 
\rangle}
f (\gamma ) \,, \]
to obtain a $ \Gamma^*$-periodic analytic function on $ \RR^2 $. 

To simplify the presentation we assume now that $ \Gamma = \ZZ^2 $. Then 
one can check \cite[\S 6]{HeSj1} the following fact:
\begin{equation}
\label{eq:op}
 \Op ( \widehat{f\; \#_B \;  g } ) = \Op (\widehat f ) \circ \Op (\widehat g ) \,, 
\end{equation}
where $ \Op $ denotes the semi-classical Weyl quantization of a function on 
$ \RR^2 $: 
\[  a ( x , \xi ) \; \longmapsto \; a^w ( x , h D) \; : \; L^2 ( \RR ) 
\; \longrightarrow \; L^2 ( \RR ) \,, \]
provided that $ a $ and all of its derivatives are bounded (see \cite{DiSj}). 
In view of \eqref{eq:mag} and \eqref{eq:op} 
it is not surprising that the invertibility of $ 
{\mathcal M}_B ( f ) $ in the algebra of magnetic matrices
is equivalent to the invertibility of $ \Op ( \widehat f ) $ in the algebra of
of pseudodifferential operators. This leads to 
\begin{thm}
Suppose that the first spectral band of a 
Schr\"odinger operator with a $\ZZ^2$-periodic smooth  potential 
is separated from other bands, with $ \theta \mapsto E(\theta ) $,
the $(2\pi \ZZ)^2$-periodic first Floquet eigenvalue. Suppose that $ P_B $
is the corresponding magnetic Sch\"odinger operator with $ B = h dx_1 \wedge 
dx_2 $. Then there exists $ ( 2 \pi \ZZ)^2 $-periodic (in $\theta$) analytic 
function, $ E = 
E ( \theta , z, h ) $, 
such that for $ z $ in a neighbourhood of the first band, and  $h $ small
\begin{gather*}
 z \in \sigma ( P_B ) \ \Longleftrightarrow  \ 
0 \in \sigma ( \Op ( E ( \bullet, z , h )) ) \,,\\
E ( \theta , z , h ) \sim E( \theta ) - z + h  E_1 (\theta ) + h^2 E_2 ( \theta, z ) 
+ \cdots \,.
\end{gather*}
\end{thm}
For the formulation for a general lattice and any dimension (in particular 
dimension three) we refer to \cite{HeSj1} where one can also find the
discussion of the coefficients in the expansion of $ E ( \theta, z , h ) $.
Considering the spectrum of the leading term, $ E( x , h D_x ) $, already shows 
how dramatic is the introduction of the 
magnetic field from the spectral point of view -- see \cite{HeSj2} and
the references given there.

\subsection{High frequency scattering by a convex obstacle}
\label{aehf}

In this section we will outline the construction of a Grushin 
problem which reduces an exterior resonance problem to 
a problem on the surface of the obstacle. It was used in \cite{SjZw0}
to describe the asymptotic distribution of resonances in scattering
by a convex obstacle satisfying a natural (at least from the point of
view of our Grushin problem) curvature pinching conditions.

The study of resonances, or scattering poles, for convex bodies 
has a very long tradition going back to Watson's 1918 work on 
electromagnetic scattering by the earth. He was motivated by 
the description of the field in the deep shadow. It provided impetus for
the work on the distribution of zeros of Hankel functions which
are the resonances for the case of the sphere. For general convex
obstacles the distribution of resonances was studied, among others, by
Buslaev, Fock, Babich-Grigoreva, Bardos-Lebeau-Rauch, and 
Harg\'e-Lebeau. We refer to \cite{SjZw0}
for pointers to the literature on the subject.

The problem can be described as follows. Let $ {\mathcal O} \subset \RR^n $ 
be a stricly convex compact set with a $ \CI $ boundary. We consider the
Dirichlet (or Neumann) Laplacian on $ \RR^n \setminus {\mathcal O} $, $ 
- \Delta_{\RR^n \setminus {\mathcal O} } $, and its resolvent,
\[ R_{\mathcal O} ( \lambda ) \defeq ( -  \Delta_{\RR^n \setminus {\mathcal O} }
- \lambda^2 )^{-1} \; : \; L^2 ( \RR^n \setminus {\mathcal O} ) 
\; \longrightarrow \; H^2 (  \RR^n \setminus {\mathcal O} )  
\cap H^1_0 ( \RR^n \setminus {\mathcal O} )  \,, \ \ \Im \lambda > 0 \,.\]
When we allow $  R_{\mathcal O} ( \lambda ) $ to act on a smaller space with
values in a larger space, it becomes meromorphic in $ \lambda $:
\[ R_{\mathcal O} ( \lambda )
\; : \; L^2_{\rm{comp}} ( \RR^n \setminus {\mathcal O} ) 
\; \longrightarrow \; H^2_{\rm{loc}} (  \RR^n \setminus {\mathcal O} )  
\cap H^1_{0,\rm{loc}} ( \RR^n \setminus {\mathcal O} )  \,, \ \
\lambda \in \left\{ \begin{array}{ll}  \CC & \text{ when $ n $ is odd} \\
\Lambda & \text{ when $ n $ is even} \end{array}\right.\, 
\]
where $ \Lambda  $ is the logarithmic plane. The poles of this meromorphic
family of operators are called {\em resonances} or {\em scattering poles}.
They constitute a natural replacement of discrete spectral data for 
problems on non-compact domains -- see \cite{Z-RPG} for an introduction and
references.

The first step of the argument is a deformation of $ \RR^n \setminus
{\mathcal O} $ to a totally real submanifold, $ \Gamma $, with boundary 
$ \partial \Gamma =\partial {\mathcal O} $ in $ \C^n $. The Laplacian $ - 
\Delta_{\RR^n \setminus {\mathcal O} }$ 
on $ \RR^n \setminus {\mathcal O} $  can be considered as a restriction of the
holomorphic Laplacian on $ \C^n $ and it in turn restricts to an operator
on $ \Gamma $, $ - \Delta_\Gamma $. When $ \Gamma $ is equal to $ e^{i 
\theta_0 } \RR^n $ near infinity then the resonances of $ -\Delta_{ 
\RR^n \setminus {\mathcal O} }$  coincide with the complex eigenvalues 
of $ \Delta_\Gamma$ in a conic neighbourhood of $ \RR $. That is the 
essence of the well known complex scaling method adapted to this setting.

Normal geodesic coordinates are obtained by taking 
$ x' $ as coordinates on $ \partial {\mathcal O} $ and 
$ x_n $ as the distance to $ \partial {\mathcal O} $.
In these coordinates the Laplacian near the boundary is approximated by 
\begin{equation}
\label{eq:o.1} D_{x_n }^2 - 2 x_n Q ( x', D_{x'} ) +  R ( x' , D_{x'} ) 
\end{equation}
where $ R $ is the induced Laplacian on the boundary and the principal
symbol of $ Q $ is the second fundamental form of the boundary.   
The complex deformation near the boundary can be obtained by rotating $ x_n $ 
in the complex plane: $ x_n \mapsto e^{i \theta } x_n $ which changes
\eqref{eq:o.1} to
\begin{equation}
\label{eq:o.2}
e^{-2i \theta} D_{x_n} ^2 - 2 e^{i \theta} x_n Q ( x' , D_{x'} ) + 
R ( x', D_{x' } ) \,.
\end{equation}
The natural choice of $ \theta $ comes from the homogeneity of the 
equation: $ \theta = \pi /3 $.

It is also natural to work in the semi-classical setting, that is,
to consider resonances of $ - h^2 \Delta_{\RR^n \setminus {\mathcal O} } $
near a fixed point, say $ 1 $. Letting $ h \rightarrow 0 $ gives then
asymptotic information about resonances of $ - \Delta_{\RR^n \setminus 
{\mathcal O} } $.

Hence we are lead to an operator which near the boundary is approximated
by 
\begin{equation}
\label{eq:o.3}
P_0 ( h ) = e^{ - 2 \pi i /3 } ( ( hD_{x_n})^2 + 2 x_n Q( x',
hD_{x'}  ) ) + R ( x' , h D_{x'}  ) \,,
\end{equation}
and we are interested in its eigenvalues close to $1$.
Let us consider the principal symbol of
\eqref{eq:o.3} in the tangential variables. That gives
\[ p_0 ( h ) = e^{ - 2 \pi i /3 } ( ( hD_{x_n})^2 + 2 x_n Q( x',
\xi' ) ) + R ( x' , \xi' ) \,.\]
We are interested in the invertibility of $ P_0( h ) - \zeta $ for
$ \zeta $ close to $ 1 $ and that should be related to invertibility
of the operator valued symbol $ p_0 ( h ) - \zeta $.
We rewrite it as
\begin{gather}
\label{eq:o.homo}
\begin{gathered}
p_0 ( h ) - \zeta = h^{\frac23} \left( e^{- 2 \pi i / 3 }
( D_t^2 +  t \mu )+ \lambda - z \right)\,, \\
t = h^{- \frac23} x_n \,, \ \lambda = h^{-\frac23} ( R ( x' , \xi') - 1 )
\,, \ z = h^{-\frac23} ( \zeta - 1) \,, \ \mu = 2 Q ( x', \xi' )\,,
\end{gathered}
\end{gather}
that is, we rescale the variables using the natural homogeneity of
$ p_0 ( h ) - \zeta $.  
On the symbolic level the operator \eqref{eq:o.3} can be analyzed rather 
easily. We can describe $ ( p_0 ( h ) - \zeta )^{-1} $ using the Airy 
function:
\[ ( D_t^2 + t ) Ai ( t ) =  0 \,, \ \ Ai ( - \zeta_j ) = 0 \,, \ \ 
Ai \in L^2( [0, \infty ) ) \,. \]

Thus we consider
\begin{equation} 
\label{eq:4.1} 
P _\lambda - z = e^{- 2\pi i / 3 } ( D_t^2 + \mu t ) + \lambda - z \,,  
 \ \  \lambda \in \RR \,, \ \  1/C \leq \mu \leq C \,, \ \
|\Im z | < C_1 \,, 
\end{equation} 
where $  C_1 $ will remain large but fixed. 
To simplify the notation we shall now put $ \mu =  1  $ (all the 
estimates will clearly be uniform 
with respect to $ \mu $ with all derivatives).
 
Let $ 0 > - \zeta_1 > - \zeta_2 > \cdots > - \zeta_k > \cdots $ 
be the zeros of the Airy function and let $ e_j ( t) = c_j  
Ai ( t -  \zeta_j )$ be the normalized eigenfunctions of  
\[ \left\{ \begin{array}{l}  
( D_t^2 + t ) e_j (t) = \zeta_j e_j ( t) \, , \ \ t \geq 0 \\ 
e_j ( 0 ) = 0 \,. \end{array} \right. 
\] 
We recall that the eigenfuctions $ e_j $ decay rapidly since for 
$ t \rightarrow + \infty  $ we have 
\[ Ai ( t ) \sim  
( 2 \sqrt {\pi } )^{-1} t^{-\frac14} \exp ( - 2 t^{\frac32} / 3 ) \,.\] 
We now take $ N = N ( C_1 ) $ the largest number such that 
\[  |\Im e^{-i 2\pi/3}  \zeta_N | \leq C_1 \,.\]  
To set up the model Grushin problem we define 
\begin{gather} 
\label{eq:4.2} 
\begin{gathered} 
R_+^0 : L^2 ( [0, \infty ) ) \;  \longrightarrow \; \C^N  \,, \ \  
R_+^0 u ( j ) = \langle u , e_j \rangle \,, \ \  1 \leq j  \leq N \,, \\ 
R_-^0 : \C^N \rightarrow L^2 ( [ 0 , \infty ) ) \,, \ \  
R_-^0 = ( R_+^0 )^*  \,.  
\end{gathered} 
\end{gather} 
Using this we put 
\begin{gather} 
\label{eq:4.new}
\begin{gathered}
 {\mathcal P}_\lambda^0 ( z) = \left(  
\begin{array}{ll} P_\lambda - z & R_-^0  \\ 
\ \ R_+^0 & 0 \end{array} \right) \; : \; B_{\lambda}
 \times \C^N  
\longrightarrow L^2 \times \C^N \,, \\ 
B_{\lambda } = \{ u \in L^2_r \; : \; D_t^2 u \,, \ t u \in L^2\,, 
\ u ( 0 ) = 0 \}
  \,, \\
\| u \|_{B_{z, \lambda  }} =  \langle \lambda - \Re z 
\rangle \| u \|_{L^2 } + 
\| D_t^2 u \|_{L^2} + \| t u \|_{L^2 } \,.
\end{gathered} 
\end{gather} 
Since the eigenvalues of $ P_\lambda  $ are given by $  
\lambda + e^{ - 2 \pi i / 3 } \zeta_j $ and $ e_j$ are the 
corresponding eigenfunctions, we see that 
 $ {\mathcal P}_\lambda  
^0 ( z ) $ is bijective with a bounded inverse. 

As in \S \ref{aeps} our Grushin problem becomes 
``stable under perturbations''. However, because of the rescaling, the symbol
class of the inverse is very bad in the original coordinates: we lose 
$ h^{-\frac23} $ when differentiating in the direction transversal to 
the hypersurface $ R - 1 = 0 $. Overcoming that requires some second microlocal
techniques. Once that is in place the invertibility of
 $  P_0 ( h ) - (1 + h^{\frac23} z )  $ 
for $ |\Im z| \leq C $  is controlled by 
invertibility of an operator on the boundary with the principal symbol
given by  
\begin{gather}
\label{eq:o.4}
\begin{gathered}
 E^0_{-+ } \in \text{Hom} ( \C^N , \C^N ) \,, \ \
( E^0_{-+} )_{1 \leq i, j \leq N} = - ( \lambda - z +
 \mu^{ \frac23} e^{ - 2 \pi i / 3 } \zeta_j ) \delta_{ij}  \,, \\ 
\lambda = h^{-\frac23} ( R ( x', \xi') - 1 )\,, \
\mu = 2 Q ( x' , \xi' ) \,. \end{gathered} \end{gather}
Here $ N $ depends on $ C $ which controls the range of $ \Im z$.

The passage to a global operator on the boundary, $ E_{-+}( z)  $, with 
poles of $ E_{-+}(z)^{-1} $ corresponding to the rescaled resonances is
rather delicate. We use \cite[Section 6]{SjZw0}  a symbolic calculus which takes
into account lower order terms near the boundary. This results in an effective
Hamiltonian, $ E_{-+} ( z ) $, described in Theorem \ref{t:7.1}.
In a suitable sense it is close to the model operator $ E_{-+}^0$
described above. It has to be stressed that a restriction on the 
range of $ \Re z $ has to be made: for every large constant $ L $ we construct
a different $ E_{-+} (z) $ which works for $ |\Re z | \leq L $. The properties
of the leading symbol remain unchanged but the lower order terms and the 
symbolic estimates depend on $L$. 

The detailed description of the effective Hamiltonian is quite technical and
involves the second microlocal classes of pseudodifferential operators 
introduces in \cite[Section 4]{SjZw0}. Nevertheless from a computational
point of view the construction is quite straightforward relying on the
Grushin problem described above and the Taylor expansion of the coefficients
of the Laplacian (in normal geodesic coordinates) at the boundary.
\begin{thm} 
\label{t:7.1} 
Let $ W \Subset ( 0 , \infty ) $ be a fixed set. For every $ w \in 
W $ and $ z \in \C $, $ |\Re z | \ll 1/ \sqrt{\delta} $, $ |\Im z | \leq C_1 $ there 
exists $ E_{w,-+ } ( z ) $, a second microlocal pseudodifferential operator 
associated to $ \Sigma_w = \{ p \in T^*\partial {\mathcal  
O} : R( p ) = w \} $, $ N = N ( C_1) $ such that for $ 0 < h < h_0  
( \delta ) $:  

\vspace{0.25cm} 
\noindent 
(i) If the multiplicity of the pole of the meromorphic continuation  
of $ ( \Delta_{\RR^n \setminus {\mathcal O } } - \zeta )^{-1} $ is given  
by $ m_{\mathcal O} ( \zeta ) $ then   
\begin{equation} 
\label{eq:7.t.1} 
  m_{\mathcal O} ( h^{-2} ( w + h^{\frac23} z )) = \frac{1}{2 \pi i } 
\; \text{\em tr} \; 
\oint_{
|\tilde z - z | = \epsilon }
E_{w, -+ } ( \tilde z )^{-1}  
\frac{d}{ d \tilde z }  E_{w, -+ } ( \tilde z ) d \tilde z  \,, 
\ \ 0  < \epsilon \ll 1\,.  
\end{equation}

\vspace{0.25cm} 
\noindent 
(ii) If $ E_{w , -+ }^0 ( z; p ) = \sigma_{ \Sigma_w ,  
h } \left( E_{ w, -+ } ( z ) \right) ( p ;
h ) $, $ p\in T^* \partial {\mathcal O} $, $ \sigma_{\Sigma_w , h } $, the
second microlocal symbol map,
\[ E_{w, -+ }^0 ( z ; p ;h  ) = {\mathcal O}  ( 
\langle \lambda - \Re z \rangle ) \,.\] 
In addition  for $ |\lambda | \leq  1 /
( C \sqrt{\delta} ) $ we have 
\begin{equation} 
\label{eq:7.t.2} 
\| E_{w , -+ }^0 ( z; p ; h) -
\text{\em diag} 
 ( z - \lambda -  
 e^{ - 2 \pi i / 3 } \zeta_j ( p ) ) \|_{ 
{\mathcal L } ( \C^N , \C^N ) } \leq \epsilon 
 \ll 1  \,,
\end{equation} 
and 
\begin{equation} 
\label{eq:7.t.3} 
\det E_{w, -+ }^0 ( z ; p ;h)  =  0 \  
\Longleftrightarrow \ z = \lambda + e^{ - 2 \pi i / 3 } \zeta_ j  
( p ) \ \text{\em for some $ 1 \leq j \leq N $}  
\end{equation}  
where the zero is simple.
Here $\zeta_j ( p ) = \zeta_j ( 2 Q ( p ) ) ^{\frac23} $. 

\vspace{0.25cm} 
\noindent 
(iii) For $ |\lambda| \geq 1/ ( C \sqrt{\delta} )  $,  $ E^0_{w , - + } $
is invertible and 
\[E^0_{w , - + }  ( z ; p ; h  ) ^{-1} = {\mathcal O} ( \langle   
\lambda - \Re z \rangle^{-1}  )\,.\]
\end{thm} 

\begin{figure}[htbp]
\begin{center}
\includegraphics[width=4.5in]{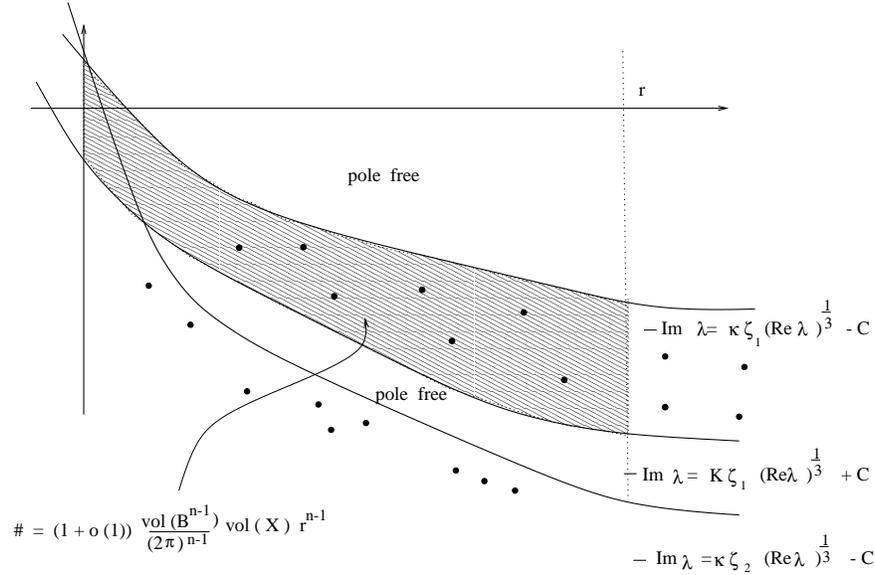}
\end{center}
\caption
{
\label{f}
The distribution of resonances for a convex obstacle satisfying the
pinched curvature assumption \eqref{eq:pinch} with $ j_0 = 1$.
}
\end{figure}

In \cite[Section 9]{SjZw0} we give  a trace formula for $ E_{-+} ( z)$. 
For that we start with the obvious observation
that the trace of the integral of 
$ E_{-+}(z)^{-1} (d/dz) E_{-+} ( z) $ against a holomorphic 
function $ f $ over a closed curve 
gives the sum of values of $ f $ at resonances enclosed by 
the curve.  The proof of the trace formula
involves a further Grushin reduction, 
a local lower modulus theorem and a good choice of contours. The gain is in 
obtaining an integral in the region where the operator $ E_{-+} ( z) $ is 
elliptic (roughly speaking in the pole free region). 
A good choice of $ f $, yields an asymptotic formula 
(see \cite[Theorem 1.2]{SjZw0}) for the number of resonances 
resonances in bands
\begin{gather*}
 \kappa \zeta_{j} (\Re \lambda)^{\frac13} - C  
 < - \Im \lambda <   K \zeta_j (\Re \lambda)^{\frac13} + C \,,  \ \ 
j \leq j_0 \\
 \kappa = 2^{-\frac13}  \cos\frac\pi6 \min_{ S\partial {\mathcal O}}  
Q^{\frac23  } \,, \ \  
K = 2^{-\frac13}  \cos\frac\pi6 \max_{ S\partial {\mathcal O}}  
Q^{\frac23  } \,, 
\end{gather*} 
where, as above, $ Q $ is the 
second fundamental form of $ \partial {\mathcal O} $ and $ S \partial 
{\mathcal O} $ the sphere bundle of $ \partial {\mathcal O}$, provided that
we have the pinched curvature condition:
\begin{equation}
\label{eq:pinch}
 \frac{ \max_{ S\partial {\mathcal O}} Q }{  \min_{ S\partial {\mathcal O}} Q }  
< \left( \frac{\zeta_{j_0+1}}{ \zeta_{j_0}} \right)^{\frac32} \,.  \end{equation}
Under this assumption the regions between the bands are resonance free -- this
is shown in Figure \ref{f} which illustrates the result.


\begin{thebibliography}{XX}


\bibitem{BaTr} D. Bau and  L.N. Trefethen, {\em Numerical Linear Algebra,}
Society for Industrial and
Applied Mathematics (SIAM), Philadelphia, PA, 1997.

\bibitem{BdM} L. Boutet de Monvel, {\em Boundary problems for 
pseudodifferential operators,} Acta Math. {\bf 126}(1971), 11-51.

\bibitem{DJ} J. Derezi\'nski and V. Jaksic,  {\em Spectral Theory of 
Pauli-Fierz Operators,} J. Funct. Anal. {\bf 180}(2001), 243--327.


\bibitem{DiSj} M. Dimassi and J. Sj\"ostrand, {\em Spectral Asymptotics in
the semi-classical limit,} Cambridge University Press, 1999.


\bibitem{KTT} R. Fletcher and T. Johnson {\em
On the stability of null-space methods for KKT systems.}
SIAM J. Matrix Anal. Appl. {\bf 18}(1997), 938--958. 


\bibitem{gk} I.C. Gohberg and M. G. Krein, {\em Introduction to the Theory 
of Linear Non-self-adjoint operators,} Translations of Mathematical 
Monographs {\bf 18}, A.M.S., Providence, 1969.

\bibitem{Gr} V.V. Grushin {\em Les probl\`emes aux limites d\'eg\'en\'er\'es 
et les op\'erateurs pseudo-diff\'erentiels.}
Actes du Congr\`es International des Math\'ematiciens 
(Nice, 1970), Tome 2, 737--743. 

\bibitem{HeSj} B. Helffer and J. Sj\"ostrand, {\em R\'esonances en limite
semi-classique. [Resonances in the semi-classical limit]}
M\'emoires de la S.M.F. {\bf 114}(3)(1986).

\bibitem{HeSj1} B. Helffer and J. Sj\"ostrand, {\em
\'Equation de Schr\"odinger avec champ magn\'etique 
et \'equation de Harper, Schrödinger operators (Soenderborg, 1988)}, 118--197, 
Lecture Notes in Phys., 345, Springer, Berlin, 1989.

\bibitem{HeSj2} B. Helffer and J. Sj{\"o}strand, {\em Semiclassical 
analysis for Harper's equation. III. Cantor structure of the spectrum.}  
M\'em. Soc. Math. France (N.S.){\bf 39}(1989), 1--124.

\bibitem{Horb1} L. H\"ormander, {\em The Analysis of Linear Partial 
Differential Operators, vol.I,II,} Springer Verlag, 1983.

\bibitem{Horb} L. H\"ormander, {\em The Analysis of Linear Partial 
Differential Operators, vol.III,IV,} Springer Verlag, 1985.

\bibitem{ISZ} A. Iantchenko, J. Sj\"ostrand, and M. Zworski,
{\em Birkhoff normal forms and semiclassical inverse problems,} 
Math. Res. Lett. {\bf 9}(2002), 337-362. 


\bibitem{Lid} V.B. Lidskii {\em 
Perturbation theory of non-conjugate operators,}
U.S.S.R. Comput. Math. and Math. Phys. {\bf 6}(1966), 73--85. 

\bibitem{MeSj} A. Melin and J. Sj\"ostrand, 
{\em Bohr-Sommerfeld quantization conditions for non-selfadjoint operators in dimension
 2}, Ast\'erisque, to appear.

\bibitem{MBO} J. Moro, J.V. Burke and M.L. Overton, 
{\em On the Lidskii-Lyusternik-Vishik Perturbation Theory for Eigenvalues with
      Arbitrary Jordan Structure,} 
      SIAM J. Matrix Anal. Appl. {\bf 18}(1997), 793-817. 

\bibitem{SjTh} J. Sj\"ostrand, {\em Operators of principal type with 
interior boundary conditions,} Acta Math. {\bf 130}(1973), 1-51.

\bibitem{Sj} J. Sj\"ostrand, {\em Pseudospectrum of differential operators,}
S\'eminaire EDP, 2002-2003, \'Ecole Polytechnique.

\bibitem{SjVo} J. Sj\"ostrand and G. Vodev, {\em Asymptotics of the number 
of Rayleigh resonances.} Math. Ann. {\bf 309}(1997), 287--306.


\bibitem{SjZw} J. Sj\"ostrand and M. Zworski, {\em Quantum monodromy and 
semi-classical trace formul{\ae}}, J. Math. Pure Appl. {\bf 81}(2002), 1-33.
See also {\em Quantum monodrom revisited,} 
{\tt math.berkeley.edu/$\sim$zworski/qmr.ps}

\bibitem{SjZw0} J. Sj\"ostrand and M. Zworski, {\em Asymptotic distribution
of resonances for convex obstacles,} Acta Math. {\bf 183}(1999), 191-253.

\bibitem{Tr} L.N. Trefethen, {\em Pseudospectra of linear operators,}
SIAM Review, {\bf 39}(1997), 383-400.

\bibitem{Zel} S. Zelditch, {\em Survey on the  Inverse Spectral Problem,}
Journal of Differential Geometry Surveys, to appear.

\bibitem{Z-RPG} M. Zworski, {\em Resonances in Physics and Geometry,}
Notices of the AMS, {\bf 46} no.3, March, 1999

\bibitem{ZwPse} M. Zworski, {\em Numerical linear algebra and 
solvability of partial differential equations,} Comm. Math. Phys.
{\bf 229}(2002), 293-307. 

\end{thebibliography}
\end{document}